
\input epsf
\magnification = 1200    
\vsize = 9 true in
\hsize = 6.8   true in 
\baselineskip = 1.3 pc    
\parskip = 1.5 pc
\parindent = 1.5 pc

\def\cov{{ <\!\!\!\!\cdot\ }}
\def\qed{{\hfill \offinterlineskip \vbox {
        \hrule
        \hbox{\vrule \phantom{i$\!$} \vrule }
        \hrule} }}

\ \vskip 1 true in 

\centerline{\bf Flag-symmetry of the poset of shuffles}
\centerline{\bf and a local action of the symmetric group}

\bigskip
\bigskip

\baselineskip = 1 pc
\halign{#&#\hfil&#&#\hfil\cr
\hskip 3.5 pc&Rodica Simion\footnote\dag{This work was carried out  
during the authors' visit to MSRI on the occasion of the
Special Year Program in Combinatorics. 
The MSRI support is gratefully acknowledged.}&\hskip 3 pc
&Richard P. Stanley\dag\footnote\ddag{Partially supported by
NSF grant DMS-9500714}\cr
&Department of Mathematics&&Department of Mathematics 2-375\cr
&The George Washington University&&Massachusetts Institute of Technology\cr
&Washington, DC 20052 &&Cambridge, MA 02139\cr
&simion@math.gwu.edu &&rstan@math.mit.edu\cr
}

\bigskip
\bigskip

\centerline{Version of February 12, 1998} 

\bigskip
\bigskip

\baselineskip = 1.1 pc

\centerline{\bf Abstract} 

{\narrower\narrower  
We show that the posets of shuffles introduced by Greene in 1988 
are flag-symmetric, and we describe a permutation action of the symmetric 
group on the maximal chains which is local and yields a representation
of the symmetric group whose character has Frobenius characteristic 
closely related to the flag symmetric function. 
A key tool is provided by a new labeling of the maximal chains
of a poset of shuffles.  This labeling and the structure of the 
orbits of maximal chains under the local action lead to combinatorial 
derivations of enumerative properties obtained originally by Greene.   
As a further consequence, a natural notion of type of shuffles 
emerges and the monoid of multiplicative functions on the poset 
of shuffles is described in terms of operations on power series.
The main results concerning the flag symmetric function and the local action
on the maximal chains of a poset of shuffles are obtained from new 
general results regarding chain labelings of posets.
}

\bigskip
\bigskip

\vfil\eject

\baselineskip = 1.3 pc

\noindent{\bf 0. Introduction}

In [St2], Stanley initiated an investigation of posets which involve 
two algebraic objects related to the order structure of the poset 
--- a certain symmetric function 
(flag symmetric function) and a certain associated representation of the 
symmetric group.  In Section 1 we give precise definitions
and summarize the results of [St2] that will be used later in this paper.
Briefly, [St2] is concerned with classes of posets whose order structure leads
to a symmetric function derived from the enumeration 
of rank-selected chains, and which turns out to be the Frobenius
characteristic of a representation of the symmetric group,
of degree equal to the number of maximal chains of the poset; 
moreover, this representation can be realized via an action of the 
symmetric group on the maximal chains of the poset, under which 
each adjacent transposition $\sigma_i = (i, i+1)$ acts on chains locally, 
that is, modifying at most the  chain element of rank $i$.   

The goal of this paper is to add a new infinite family of posets to
the examples appearing in [St2] and [St3], namely, the posets of
shuffles introduced and investigated by Greene [Gre1].  In the
process, several general results emerge.  In Section 1 we give the
necessary background on locally rank-symmetric posets affording a
local action of the symmetric group (based on [St2]).  Section 2
contains the necessary preliminaries concerning the posets of
shuffles 
(i.e., shuffles of subwords of two given words).  
In Section 3 we give a new labeling of the posets of
shuffles and establish its properties which are instrumental in the
remainder of the paper.  In Section 4 we describe a local action of
the symmetric group on the maximal chains of a poset of shuffles, such
that the Frobenius characteristic for the corresponding representation
character is (essentially) the flag symmetric function.  The desired
results regarding the posets of shuffles follow from more general
results motivated by the properties of the new labeling of these
posets.  Section 5 is devoted to the enumeration of shuffles according
to a natural notion of type.  As a consequence we describe the monoid
of multiplicative functions on the poset of shuffles in terms of
operations on power series.

As a by-product of the present investigation of the posets of
shuffles, we obtain alternative, purely combinatorial, derivations of
enumerative results obtained in [Gre1].  The present work parallels
that of [St3] regarding the lattice of noncrossing partitions, thus
adding to previously known structural analogies between the posets of
non-crossing partitions and those of shuffles.  It is hoped that this
work will facilitate the development of a systematic general theory of
the posets with a local group action concordant with the flag
symmetric function.

\bigskip
\bigskip

\noindent{\bf 1. Preliminaries}

Let $P$ be a finite poset with a minimum element $\hat 0$, 
and a maximum element $\hat 1$.  Throughout this paper, we will consider 
only such posets that are {\it ranked}, that is, there exists a function 
$\rho \colon P \to {\bf Z}$ 
such that $\rho({\hat 0}) = 0$ 
and $\rho(t') = \rho(t) + 1$ whenever $t \cov t'$ (the notation 
$t \cov t'$ means that $t$ is {\it covered} by $t'$, i.e., $t<t'$ and
there is no element $u \in P$ such that $t < u < t'$). 

Let $\rho(P)\ \colon = \rho({\hat 1}) = n$.  For $S \subseteq [n -1]$, 
where $[n-1]\ \colon = \{ 1, 2, \dots, n-1 \}$, let $\alpha_P(S)$ denote 
the number of rank-selected chains in $P$ whose elements (other than 
$\hat 0$ and $\hat 1$) have rank set equal to $S$.  Thus, 
$$
\alpha_P(S) \colon = \ \# \{ {\hat 0} < t_1 < t_2 < \cdots < t_{|S|} < {\hat 1}
\ \colon \ \{ \rho(t_1), \rho(t_2), \dots, \rho(t_{|S|}) \} = S \} .
$$
The function $\alpha_P \colon 2^{[n-1]} \to {\bf Z}$ 
is the {\it flag f-vector} of $P$.
It contains information equivalent to that of the {\it flag h-vector} 
$\beta_P$ whose values are given by 
$$
\hskip 1 true in 
\beta_P(S) \colon = \ \sum_{T \subseteq S}{(-1)^{|S - T|} \alpha_P(T)},
\ \hskip .7 true in {\rm for\ all \ } S \subseteq [n-1].
\hskip .7 true in (1)
$$
For example, writing $\alpha_P(t_1,\dots,t_k)$ for $\alpha_P(\{
t_1,\dots, t_k\})$ and similarly for $\beta_P$, 
the poset of Figure 1 has 
$\alpha(\emptyset ) = 1, \alpha(1 ) = \alpha(  2  ) =  5, 
\alpha(  1, 2  ) = 12$, and 
$\beta(\emptyset ) = 1, \beta({1 }) = \beta({  2  }) = 4, 
\beta({  1, 2  }) = 3$.  The poset of Figure 2 
has 
$\alpha (\emptyset ) = 1, \alpha ( 1  ) = \alpha (  3  ) = 2, 
\alpha (  2  ) = 3, 
\alpha (  1, 2  ) = \alpha (  1, 3  ) = \alpha (  2, 3  ) = 4,
\alpha (  1, 2, 3  ) = 6$, and 
$\beta(\emptyset ) = 1, \beta({1 }) = \beta({  3  }) = 1, 
\beta({  2  }) = 2, 
\beta({  1, 2  }) = \beta({  2, 3  }) = 0,
\beta({  1, 3  }) = 1,
\beta({  1, 2, 3  }) = 0$.

%
%

        \midinsert
        \epsfysize=3in
       \centerline{ \epsfbox{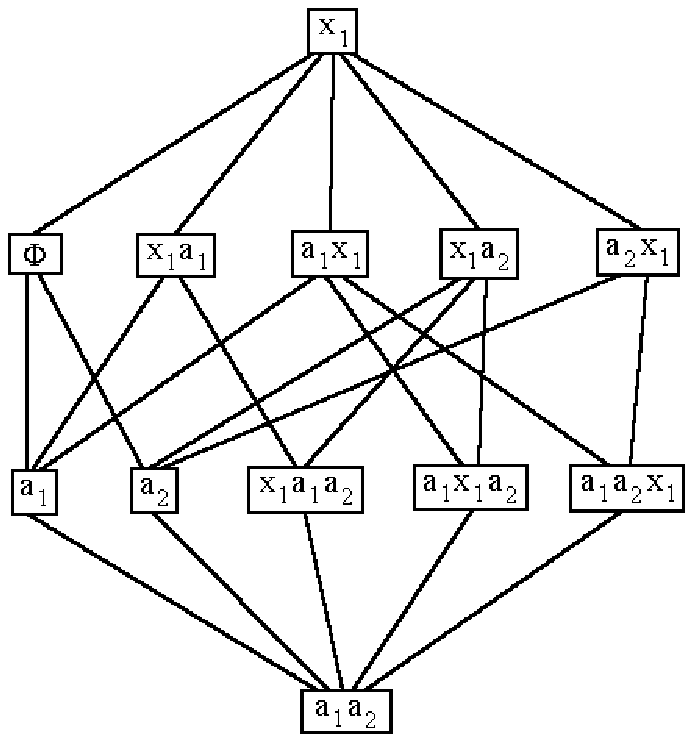} }
        \vskip -3ex
        \centerline {{\bf Figure 1.} The poset of shuffles $W_{21}$.}
        \endinsert

The flag $f$- and $h$-vectors appear in numerous contexts in  
algebraic and geometric combinatorics;  for instance, the values 
$\beta_P(S)$ have topological significance related to the order complex of 
the rank-selected subposet 
$P_S \colon = \{ {\hat 0}, {\hat 1} \} \cup 
              \{ t \in P \ \colon \ \rho(t) \in S \}$ 
(see, e.g.,  [St1, section 3.12] for additional information and references). 
   
Consider now the formal power series 
$$
F_P(x) \colon = F_P(x_1, x_2, \dots ) = 
\sum_{ {\hat 0} \leq t_1 \leq t_2 \leq \cdots \leq t_{k} < {\hat 1} }
     { x_1^{\rho(t_1)} x_2^{\rho(t_2) - \rho(t_1)} \cdots 
       x_{k+1}^{n - \rho(t_k)} }.$$
This definition was suggested for investigation by Richard Ehrenborg
[E] and is one of the central objects in [St2] and in this paper. 
Alternatively, 
$$
\hskip .4 true in
F_P(x) = \sum_{{S \subseteq [n-1]} \atop{ S = \{ s_1 < s_2 < \cdots < s_k \}} }
     { \alpha_P(S) \ \cdot 
       \sum_{1 \leq i_1 < i_2 < \cdots < i_{k+1} } 
            {x_{i_1}^{s_1} x_{i_2}^{s_2 - s_1} \cdots x_{i_{k+1}}^{n - s_k} }
      }.
\hskip 0.65 true in (2)
$$ 
It is easy to see that the series $F_P(x)$ is homogeneous of degree $n$
and that it is a {\it quasisymmetric function}, that is, 
for every sequence $n_1, n_2, \dots, n_m$ of exponents, 
the monomials $x_{i_1}^{n_1} x_{i_2}^{n_2} \cdots x_{i_m}^{n_m}$ 
and  $x_{j_1}^{n_1} x_{j_2}^{n_2} \cdots x_{j_m}^{n_m}$ 
appear with equal coefficients whenever 
$i_1 < i_2 < \cdots < i_m$ and 
$j_1 < j_2 < \cdots < j_m$.
Through a simple counting argument 
and using the relation (1), the series $F_P(x)$ can also be rewritten as 
$$
\hskip 1.65 true in
F_P(x) = \sum_{S \subseteq [n-1]} { \beta_P(S)\ L_{S,n}(x) },
\hskip 2 true in (3)
$$
where the $L_{S,n}(x)$ are Gessel's quasisymmetric functions 
$$
L_{S,n}(x) \ \colon = \sum_{ {1 \leq i_1 \leq i_2 \leq \cdots \leq i_n}
                           \atop { i_j < i_{j+1} \ {\rm if} \ j \in S} }
                    {x_{i_1} x_{i_2} \cdots x_{i_n}},  
$$
which constitute a basis for the ($2^{n-1}$-dimensional) space of 
quasisymmetric functions of degree $n$ (for more on quasisymmetric
functions  and symmetric functions we refer the interested reader to 
[Re] and [M]). 

A first question discussed in [St2] is that of conditions under which 
$F_P(x)$ is actually a symmetric function, in which case we refer 
to $F_P(x)$ as the {\it flag symmetric function of $P$}
and to $P$ as a {\it flag-symmetric poset}. 
An immediate necessary condition is that $P$ be rank-symmetric 
(i.e., $\# \{ t \in P \ \colon \ \rho(t) = r \}  = 
        \# \{ t \in P \ \colon \ \rho(t) = n-r \}  $ for every 
$0 \leq r \leq n$). 
A necessary and sufficient condition can be deduced readily from (2)
[St2, Corollary 1.2].  
Namely, for every $S \subseteq [n-1]$ the value of $\alpha_P(S)$ 
depends only on the (multi)set of differences 
$s_1 -0, s_2 - s_1, \dots, s_{k} - s_{k-1}, n-s_k$ 
and not on their ordering.  If this is the case, then the symmetric function 
$F_P(x)$ can be expressed in terms of the basis of monomial symmetric 
functions $\{ m_{\lambda}(x) \}_{\lambda \vdash n}$ as 
$$
\hskip 1.7 true in
F_P(x) = \sum_{\lambda \vdash n} 
         { \alpha_P({\lambda}) m_{\lambda}(x) },
\hskip 2 true in (4)
$$
where $\lambda \vdash n$ denotes a partition 
$\lambda = (\lambda_1 \ge \lambda_2 \ge \cdots \ge \lambda_l > 0)$ 
of $n$, 
and $\alpha_P({\lambda}) \colon = 
\alpha_P({ \{ \lambda_1, \lambda_1 + \lambda_2, \dots, 
\lambda_1 + \lambda_2 + \dots + \lambda_{l-1} \} })$. 

%

        \midinsert
        \epsfysize=3in
       \centerline{ \epsfbox{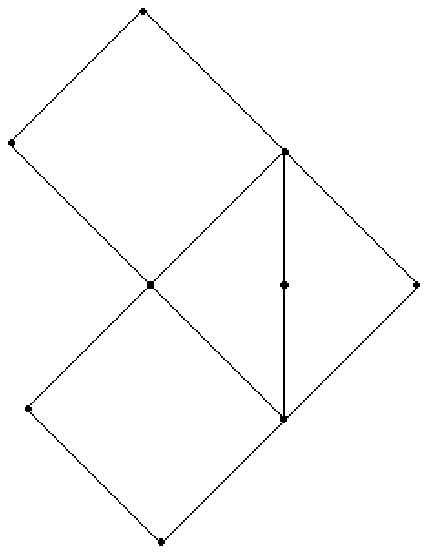} }
        \vskip -3ex
        \centerline {{\bf Figure 2.}  A flag-symmetric poset which is
     not locally  rank-symmetric.}
        \endinsert

For example, the earlier calculation of the 
flag $f$-vector of the poset of Figure~2 shows that 
it is flag-symmetric, with flag symmetric function
$F_P(x) = m_{4}(x) + 2 m_{31}(x) + 3 m_{22}(x) + 4 m_{211}(x) + 
6 m_{1111}(x)$.

The following sufficient condition for $F_P(x)$ to be 
a symmetric function introduces the class of {\it locally rank-symmetric
posets}.
This condition is not necessary for flag-symmetry, 
as shown by the poset of Figure 2, but it is necessary and sufficient for 
every interval of $P$ to be flag-symmetric. 

\noindent {\bf Proposition 1.1} [St2, Theorem 1.4] 
{\it 
Let $P$ be a ranked poset with $\hat 0$ and $\hat 1$.
If $P$ is locally rank-symmetric, i.e., if every interval in $P$ 
is a rank-symmetric (sub)poset, then $F_P(x)$ is a symmetric function. 
}

Locally rank-symmetric posets turn out to be a rich source of examples
yielding flag symmetric functions.  The examples of flag-symmetric posets
provided in [St2] include products of chains (shown to be the only 
flag-symmetric distributive lattices, and identical to the class 
of locally rank-symmetric distributive lattices), and Hall lattices (a
``$q$-analogue'' of a product of chains),
as well as a discussion of some other classes of posets.
If $F_P(x)$ is a symmetric function,
homogeneous of degree $n$,
and if it turns out to be Schur positive
(i.e., its expression in terms of the Schur functions basis has 
nonnegative coefficients only), then it follows from the general theory of 
representations and symmetric functions that it is the Frobenius
characteristic 
$$
\hskip 1.85 true in
 {\rm ch}(\psi) \  \colon = \sum_{\lambda \vdash n} 
            { { {\psi(\lambda)} \over {z_{\lambda}} } p_{\lambda}(x) }
\hskip 1.95 true in (5)
$$
of a character $\psi$ of the symmetric group $S_n$. 
In the preceding display line, $\lambda$ runs over all partitions of $n$, 
$\psi(\lambda)$ is the value of $\psi$ on the conjugacy class of type 
$\lambda$, $z_{\lambda} = 1/(\lambda_1 \lambda_2 \cdots m_1! m_2! \cdots )$
with $m_i$ being the multiplicity of $i$ as a part of $\lambda$, 
and $p_{\lambda}(x)$ 
is the power symmetric function indexed by $\lambda$ 
(that is, $p_{\lambda}(x) = p_{\lambda_1}(x) p_{\lambda_2}(x) \cdots$ 
with $p_j (x) \colon = x_1^j + x_2^j + \cdots$). 
It is known that when the Frobenius characteristic 
of a character $\psi$ of $S_n$ is expanded in terms 
of Schur functions $\{ s_{\lambda}(x) \}_{\lambda \vdash n}$, 
then the coefficient of $s_{\lambda}(x)$ is 
the multiplicity with which the irreducible character $\chi^{\lambda}$ 
of $S_n$ occurs in $\psi$.  
Thus, $F_P(x)$ describes a representation of $S_n$,
whose degree $\psi(1^n)$ can be recovered as  
the coefficient of 
$m_{1^n}$ in $F_P(x)$. 
In view of (4), 
the degree of $\psi$ is $\alpha_P({1^n})$, the number of maximal chains in $P$. 

The preceding discussion suggests seeking a natural action of $S_n$ on
the complex vector space ${\bf C}{\cal M}(P)$  with the set 
${\cal M}(P)$ of maximal chains in $P$ as a basis, giving  rise 
to a representation of $S_n$ with character $\psi$ as in (5). 
Of particular interest would be a {\it local action} with this property 
(defined in [St2] and motivated by the notion of local stationary algebra
appearing  in [V]);  that is, an action such that 
for every adjacent transposition $\sigma_i = (i, i+1)$ and every maximal
chain $m$ of $P$ we have 
$$
\sigma_i(m) = \sum_{m' \in {\bf C}{\cal M}(P) }
                  { c_{m m'} m' },
$$
with nonzero coefficient $c_{m m'}$ only if $m'$ differs from $m$
at most in the element of rank $i$. 
Following [St2], we call such an action {\it good}. 
Good actions of the symmetric group are discussed in [St2]
in the case of  posets whose rank-two intervals
are isomorphic to $C_3$ or $C_2 \times C_2$ (where $C_i$ denotes an
$i$-element chain), 
and for posets whose rank-three intervals are isomorphic to 
$C_4$ or $C_3 \times C_2$ or $C_2 \times C_2 \times C_2$. These
results are based on work of David Grabiner [Gra].
Another illustration in [St2] gives a  
local action of the Hecke algebra of $S_n$
on ${\bf C}{\cal M}(B_n(q))$, where $B_n(q)$ denotes the lattice of subspaces 
of an $n$-dimensional vector space over $GF_q$.  
In [St3] a good action is exhibited for the lattice of noncrossing
partitions.
To these classes of examples this paper adds the posets of shuffles.

We note that related results were recently obtained by Patricia Hersh
[He] (generalizing the local $S_n$-action on noncrossing partitions),
and Jonathan Farley and Stefan Schmidt [FaSc] (generalizing the work
of Grabiner [Gra]).


\bigskip
\bigskip

\noindent{\bf 2. Flag-symmetry of the posets of shuffles}

Let ${\cal A} = \{ a_1, a_2, \dots, a_M \}$ and 
 ${\cal X} = \{ x_1, x_2, \dots, x_N \}$
be two (finite) disjoint sets which we will call the lower 
and the upper alphabets, respectively.  Consider the collection of 
{\it shuffles} over ${\cal A}$ and ${\cal X}$, that is, words 
$w = w_1 w_2 \cdots w_k$ with distinct letters from ${\cal A} \cup {\cal X}$ 
satisfying the {\it shuffle property:}
the subset of letters belonging to each alphabet appears in 
increasing order of the letter-subscripts in the appropriate alphabet. 
For instance, if $M=4$ and $N = 3$, then $w = x_2 a_1 a_3 x_3$ is 
a shuffle word, but $w = a_1 x_2 a_2 a_3 x_1$ is not a shuffle word. 
Note that the empty word $\emptyset$ is a shuffle word.

The {\it poset of shuffles} $W_{MN}$ 
consists of the shuffle words over alphabets $\cal A$ and $\cal X$ 
with $\# {\cal A} =M$ and $\# {\cal X} = N$ with the order relation
given by $w \cov w'$ iff $w'$ is obtained from $w$ either 
by deleting a letter belonging to $\cal A$ or by 
inserting (in an allowable position) a letter belonging to $\cal X$. 
In particular, ${\hat 0} = a_1 a_2 \cdots a_M$ 
and ${\hat 1} = x_1 x_2 \cdots x_N$. 
Figure 1 shows the Hasse diagram of $W_{2 1}$. 
Clearly, $W_{0 N}$ and $W_{M 0}$ are isomorphic to 
the boolean lattices of rank $N$ and $M$, respectively. 
We will write $\{ w \}$ for the set of letters of a shuffle word $w$.

Greene [Gre1] investigated the posets of shuffles, whose definition was 
motivated by an idealized model considered in mathematical biology.
The results established in [Gre1] include structural properties of $W_{M N}$ 
(e.g., $W_{M N}$ is a ranked poset; it admits a decomposition into 
symmetrically embedded boolean lattices and, hence, 
a symmetric chain decomposition;  $W_{M N}$ is an EL-shellable poset),
as well as expressions for key invariants of $W_{M N}$ 
(the zeta polynomial, the number of maximal chains, 
 the M\"obius function, the rank generating function, 
 the characteristic polynomial). 
Two of the formulas in [Gre1] will arise later in this paper. 

\noindent
{\bf Proposition 2.1.} [Gre1, Theorem 3.4] 
{\it
The number of maximal chains in $W_{M N}$ is given by 
$$
\hskip 1.17 true in
C_{MN} = (M+N)! \ \sum_{k \ge 0} 
         { {M \choose k} {N \choose k} { 1 \over {2^k} } }.
\hskip 1.5 true in (6)
$$
The M\"obius function of $W_{M N}$ is 
$$
\hskip 1 true in
\mu_{MN} = \mu_{W_{M N}}({\hat 0}, {\hat 1}) = (-1)^{M+N} {{M+N}\choose M}.
\hskip 1.3 true in (7)
$$
} 

We now turn to the interval structure of the posets of shuffles.

\noindent{\bf Lemma 2.2.} 
{\it 
Every interval in a poset of shuffles is isomorphic to a product 
of posets of shuffles.
}

\noindent{\bf Proof.}
   Let $[u,w]$ be an arbitrary interval in $W_{M N}$, 
   and write $u = u_1 u_2 \cdots u_r$, $w = w_1 w_2 \cdots w_s$.
   Let $u_{i_1} u_{i_2} \cdots u_{i_t}$ 
   and $w_{j_1} w_{j_2} \cdots w_{j_t}$ be the subwords of $u$ and $w$,
   respectively, formed by the letters common to the two words. 
   Because $u < w$,  the shuffle property implies $u_{i_p} = w_{j_p}$ 
   for each $p = 1, 2, \dots, t$.  Moreover, the remaining letters 
   of $u$ belong to the alphabet ${\cal A}$ and the remaining letters
   of $w$ belong to the alphabet ${\cal X}$. 
   Therefore the interval $[u,w]$ is isomorphic to the product of the 
   posets of shuffles $W_{i_p - i_{p-1} -1, j_p - j_{p-1} -1}$ 
   for $p = 1, 2, \dots , t+1$, where we set $i_0 = j_0 = 0$, 
   $i_{t+1} = r+1$ and $j_{t+1} = s+1$. 
\qed

For example, if $u = a_2 x_3 a_4 a_5 a_{10} x_6 x_8$ and  
   $v = x_1 x_2 x_3 x_5 a_{10} x_6 x_8 x_{10} x_{11}$ 
   in $W_{10, 15}$, 
   then  $r=7, s=9$ and there are $t = 4$ letters common to the two words.
   These form the word 
   $x_3 a_{10} x_6 x_8 = u_2 u_5 u_6 u_7 =  v_3 v_5 v_6 v_7$
   so we have $[u, v] 
   \simeq W_{1 2} \times W_{2 1} \times W_{0 0} \times W_{0 0} \times W_{0 2}  
   \simeq W_{1 2} \times W_{2 1} \times W_{0 2}$. 

\def\caniso{{  {\simeq}_{\bf c} }}  

\noindent{\bf Remark 2.3.}
 Of course, factors of the form 
   $W_{0 0}$ are singleton posets and can be discarded from the product,
and $W_{i 0} \simeq W_{0 i} \simeq B_i$, the boolean lattice with $i$ atoms.
Using the notation from the proof of Lemma 2.2, we will write 
$[u, w] {\caniso}  \prod_{p}{ W_{i_p - i_{p-1} -1, j_p - j_{p-1} -1}}$,
the {\it canonical isomorphism type} of the interval $[u,w]$.
The notion of canonical isomorphism type of an interval will 
be used in Section 5. 
  
\noindent{\bf Proposition 2.4.}
{\it 
For every $M,N \ge 0$, the poset of shuffles $W_{M N}$ 
is locally rank-symmetric.
}

\noindent{\bf Proof.}
   Since the posets of shuffles are rank-symmetric [Gre1, Corollary 4.9]
   and since the product of posets preserves rank-symmetry, 
   Lemma 2.2 implies that  that every interval in $W_{M N}$  
   is rank-symmetric.
   \qed

 It therefore follows from Proposition~1.1 [St2, Theorem 1.4] that
each  poset of shuffles $W_{M N}$ has a flag symmetric
 function $F_{M N} = F_{M N}(x) \colon = F_{W_{M N}}(x)$.
An explicit  expression for $F_{M N}$ can be 
obtained by extending Greene's notion of ${\cal A}$- 
 and ${\cal X}$-maximal chains in $W_{M N}$ [Gre1] to 
${\cal A}$- and ${\cal X}$-maximal chains in rank-selected subposets 
of $W_{M N}$.  An argument similar to Greene's yields a recurrence 
relation for the numbers $\alpha_{W_{M N}}(\lambda)$,
which in turn implies that 
$$
\hskip .7 true in
F_{M N} = F_{M-1, N}\, e_1 + F_{M, N-1}\, e_1 - 
F_{M-1, N-1}\, e_2 - F_{M-1, N-1}\, p_2,
\hskip 1.1 true in (8)
$$
leading to 
$$
\hskip 1.4 true in
\sum_{M,N \ge 0}{F_{M N} u^M v^N } = 
   { 1 \over {(1- u  e_1)(1 - v  e_1) - u v  e_2}},
\hskip 1.5 true in (9)
$$
   where $e_j = \sum_{1 \leq i_1 < i_2 < \cdots < i_j}
                        {x_{i_1} x_{i_2} \cdots x_{i_j}},$
   the $j$th elementary symmetric function in variables $x_1, x_2,
   \dots$, and $p_2=\sum x_i^2$. Consequently, 
$$
\hskip 1.8 true in
F_{M N} = \sum_{k \ge 0} { 
            { M \choose k} {N \choose k} e_2^k\, e_1^{M+N-2k}}.
\hskip 1.95 true in (10)
$$
For example, the calculation of the flag $f$-vector of $W_{2 1}$ 
done in Section 1 gives $F_{2 1}(x) = m_3(x) + 5 m_{2 1}(x) + 12 m_{1 1 1}(x)$.
Since $e^3_1(x) = m_3(x) + 3 m_{2 1}(x) + 6 m_{1 1 1}(x)$ 
and $e_2(x) e_1(x) = m_{2 1}(x) + 3 m_{1 1 1}(x)$, 
we have $F_{2 1}(x) = e^3_1(x) + 2 e_2(x) e_1(x)$.
 
We omit the details of this argument.  Instead, we will 
obtain the expression (10) for the flag-symmetric function of 
a poset of shuffles as a consequence (Corollary 4.5) of a general result 
(Theorem 4.4) concerning chain labelings of flag-symmetric posets.

%
%
%
%
%
%
%
%
%


\bigskip
\bigskip


\bigskip
\bigskip

\noindent{\bf 3.  A labeling for posets of shuffles}

To describe an action of $S_{M+N}$ on the maximal chains of $W_{M N}$, 
it would be natural to resort to a labeling of the chains 
and have the symmetric group act on the chains by acting on 
their label sequences simply by permuting coordinates.  The poset
$W_{M N}$ is already known  
to be EL-shellable [Gre1], through the labeling of each covering relation 
 $u \cov w$ by the unique letter in the symmetric difference 
of the sets of letters $\{ u \}$ and $\{ w \}$, 
and with the ordering $a_1 < a_2 < \cdots < a_M < x_1 < x_2 \cdots < x_N$
for the labels.  Under this labeling 
each maximal chain is labeled by a permutation in $S_{M+N}$. 
However, this does not serve well the goal 
of describing an $S_{M+N}$ action on the maximal chains. 
A similar situation occurred in [St3], where the standard EL-labeling
of the noncrossing partition lattice was not suitable for 
describing a local action of the symmetric group on the maximal chains,
and a new EL-labeling was produced for this purpose. 
Here too, we will define a new labeling for a poset of 
shuffles which lends itself naturally to the 
description of the desired local action of $S_{M+N}$. 

By a {\it labeling} we mean a map $\Lambda \colon {\cal M}(P) \to
L^n$, written
 $$ \Lambda(c) = (\Lambda_1(c),\Lambda_2(c),\dots,\Lambda_n(c)), $$
where $n$ is the length of the maximal chains of $P$, and $L$ is a
totally ordered set.  The labeling of interest in the present paper is
a {\it $C$-labeling}, that is, for every maximal chain $c = (
{\hat{0}} = w^0 \cov w^1 \cov \cdots \cov w^n = {\hat 1})$ and every
$r \in [n]$, the label $\Lambda_{r}(c)$ depends only on the initial
subchain $( {\hat{0}} = w^0 \cov w^1 \cov \cdots \cov w^{r})$.  If the
label $\Lambda_r(c)$ depends only on the covering $w^{r-1} \cov w^r$,
and not on the maximal chain $c$ itself, then $\Lambda$
is an {\it $E$-labeling}.

Three properties of labelings will play a role in this paper:  
the $R^*$-, $R$-, and $S$-labeling properties. 
A $C$-labeling $\Lambda$ is an {\it $R^*$-labeling} 
if every chain of the form 
$$
( {\hat{0}} = w^0 \cov w^1 \cov \cdots \cov w^r < u)
$$
has a unique completion by covering relations 
$$
( {\hat{0}} = w^0 \cov w^1 \cov \cdots \cov w^r
\cov w^{r+1} \cov \cdots \cov w^s= u)
$$
such that 
 $$ \Lambda_{r+1}(c)<\Lambda_{r+2}(c)<\cdots<\Lambda_s(c), $$
where $c$ is any maximal chain beginning $\hat{0}=w^0\cov w^1\cov
\cdots\cov w^s$. (By the definition of $C$-labeling, the remaining
elements of $c$ do not affect the labels $\Lambda_i(c)$ for $1\leq
i\leq s$.)
In the same setting as for an $R^*$-labeling, 
the requirement for an {\it $R$-labeling} 
is the existence of a unique {\it weakly} increasing completion of the
chain:
 $$ \Lambda_{r+1}(c)\leq\Lambda_{r+2}(c)\leq\cdots\leq\Lambda_s(c). $$
A labeling $\Lambda$ 
of the maximal chains of a poset is an {\it $S$-labeling} 
if it is one-to-one and if for every maximal chain 
$c =  ( {\hat{0}} = w^0\cov w^1 \cov \cdots \cov w^n= {\hat{1}} )$
and for every rank $i \in [n-1]$ such that 
$\Lambda_i(c)\neq\Lambda_{i+1}(c)$,
there is a unique chain 
$c' =  ( {\hat{0}} = w^0 \cov w^1 \cov \cdots 
\cov w^{i-1} \cov t^i \cov w^{i+1} \cov \cdots 
\cov w^n = {\hat{1}} )$
differing from $c$ only at rank $i$, with the following property: the
label sequence  
$\Lambda(c')$ differs from $\Lambda(c)$ 
only in that $\Lambda_i(c')=\Lambda_{i+1}(c)$ and $\Lambda_{i+1}(c') =
\Lambda_i(c)$.

We now turn to the definition of a labeling $\Lambda$ 
for the poset of shuffles, and then show that it has the properties 
$R^*$ and $S$.  In the next section we will see the implications 
of an $RS$- or  $R^*S$-labeling with regard to a local action on the poset. 

%
%
%

        \midinsert
        \epsfysize=3in
       \centerline{ \epsfbox{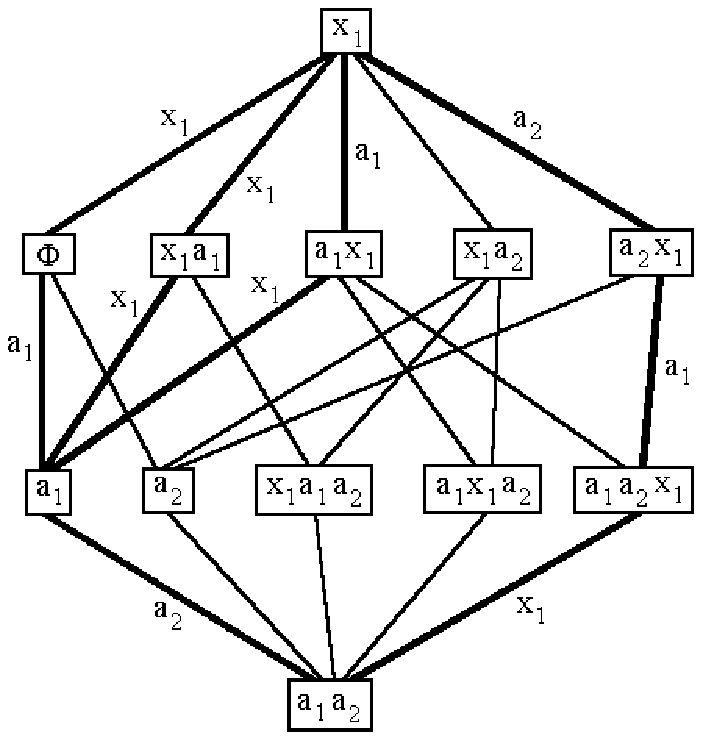} }
        \vskip -1ex
        \centerline {{\bf Figure 3.} The labeling $\Lambda$}
    \centerline{on four of the maximal chains of $W_{2 1}$.} 
        \endinsert

To each maximal chain 
$c \ = \ 
( {\hat{0}} = w^0 \cov w^1 \cov \cdots \cov w^{M+N}= {\hat{1}} )$
in $W_{M N}$ we give a label sequence
$$
\Lambda(c) = ( \Lambda_1(c),\Lambda_2(c),\dots,\Lambda_{M+N}(c)), 
$$
by assigning a label from ${\cal A} \cup {\cal X}$ 
to each covering relation on $c$.  In defining $\Lambda$  
we distinguish three types 
of covering relations, ($x$), ($xa$), and ($a$), as follows:
\indent\item{($x$)} $w^i \cov w^{i+1}$ with $w^{i+1}$ 
obtained from $w^i$ by inserting a letter $x_k \in {\cal X}$
in a position consistent with the shuffle property;
then we set $\Lambda_{i+1}(c) = x_k$. 
\indent\item{($xa$)} $w^i \cov w^{i+1}$ with 
$w^i$ of the form $w^i = u x_k a_m v$
and $w^{i+1} = u x_k v$, where this is the first deletion 
along $c$, starting from $\hat 0$, 
of a letter (necessarily belonging to $\cal A$) 
located immediately after $x_k$; 
then we set $\Lambda_{i+1}(c) = x_k$. 
\indent\item{($a$)} $w^i \cov w^{i+1}$ with 
$w^{i+1}$ obtained from $w^i$ by deleting a letter $a_j \in {\cal A}$,
and this deletion is not of type ($xa$);
then we set $\Lambda_{i+1}(c)=a_j$.
 
 Figure 3 shows the labeling $\Lambda$ of four of the maximal 
 chains in $W_{2 1}$.

\noindent{\bf Lemma 3.1.} 
{\it
The labeling $\Lambda$ is injective on the maximal chains of $W_{M N}$, 
for all $M, N$.
Its range consists of the (multi)permutations of 
all multisets of the form $A \cup 2 X \cup ({\cal X} - X)$, 
where $A \subseteq {\cal A}$, $X \subseteq {\cal X}$, 
$|A| + |X| = M$, and $2 X$ denotes the multiset consisting of two 
copies of each element of $X$. 
}

{\bf Proof.} 
{From} the definition of $\Lambda$ it is clear that 
all letters in $\cal X$ appear in the label sequence of any maximal chain $c$ 
and that for every $a_j \in {\cal A}$ which does not appear  
in the label sequence, there is an $x_k \in {\cal X}$
which appears twice.
Thus, the label sequence of every maximal chain $c$ is of the claimed form.

Conversely, we claim that given a multipermutation $\sigma$ of  
$A \cup 2 X \cup ({\cal X} - X)$ for some $A$ and $X$ as in the 
statement of the Lemma, there is a unique maximal chain
in $W_{M N}$ having label sequence $\Lambda(c) = \sigma$. 
Indeed, first note that if $A = {\cal A}$ and $X = \emptyset$ 
(that is, $\sigma$ is a permutation of ${\cal A} \cup {\cal X}$),
then only coverings of type ($a$) and ($x$) are possible.
Thus, starting from ${\hat 0} = w^0$, 
$\sigma$ dictates a sequence of deletions of letters from $\cal A$ 
and insertions of letters from $\cal X$, each insertion being made 
in the rightmost possible position.  This determines a unique 
maximal chain $c$ as desired. 
For example, for $W_{2 3}$, the permutation 
$\sigma = a_2 x_3 x_1 a_1 x_2$ 
determines the chain 
$c = ({\hat 0} = a_1 a_2 \cov a_1 \cov a_1 x_3 \cov a_1 x_1 x_3 \cov 
x_1 x_3 \cov x_1 x_2 x_3 = {\hat 1})$. 

Next, 
suppose that $ {\cal A} - A = \{ a_{i_1}, a_{i_2}, \dots, a_{i_k} \}$ 
and $X = \{ x_{j_1}, x_{j_2}, \dots, x_{j_k} \}$, 
for  some $1 \leq k \leq {\rm min} \{ M, N \}$,
where  ${i_1} < {i_2} < \cdots < {i_k}$  
and ${j_1} < {j_2} < \cdots < {j_k}$.
Observe that the shuffle condition implies that if 
the pairs $x_m, a_p$ and $x_n, a_q$ are involved in coverings 
of type ($xa$), then $m \neq n$ and $p \neq q$, 
and $m < n$ if and only if $p < q$. 
Therefore, in the multipermutation $\sigma$, 
the second occurrence of $x_{j_r}$ must correspond to a covering 
of type ($xa$) involving the pair of letters $x_{j_r}, a_{i_r}$,
for each $r = 1, \dots, k$. 
The first occurrence of $x_{j_r}$ in $\sigma$ is forced to 
correspond to the insertion of $x_{j_r}$ immediately in front of 
$a_{i_r}$, and for each $x_t \not\in X$, 
its unique occurrence in $\sigma$ forces the insertion of $x_t$
in the rightmost position possible to the left of $a_{i_r}$
and/or $x_{j_r}$, if $j_r = {\rm min} \{ s > t \ \colon \ x_s \in X \}$
(if this set is empty, then $x_t$ is inserted in the rightmost 
position possible). 
As in the preceding case, a unique maximal chain $c$ is 
determined by $\sigma$. 
For example, for $W_{4 5}$, let 
$\sigma = x_3 x_5 a_2 x_4 x_2 x_1 x_2 x_4 a_4$.
The coverings of type ($xa$) must involve the pairs 
$x_2, a_1$ and $x_4, a_3$. 
{From} $\sigma$ we reconstruct the chain 
\hfil\break
${\hat 0} = a_1 a_2 a_3 a_4 \cov a_1 a_2 x_3 a_3 a_4 \cov 
a_1 a_2 x_3 a_3 a_4 x_5 \cov a_1 x_3 a_3 a_4 x_5 \cov 
a_1 x_3 x_4 a_3 a_4 x_5 \cov 
$
\hfil\break
$x_2 a_1 x_3 x_4 a_3 a_4 x_5 \cov 
x_1 x_2 a_1 x_3 x_4 a_3 a_4 x_5 \cov 
x_1 x_2 x_3 x_4 a_3 a_4 x_5 \cov 
x_1 x_2 x_3 x_4 a_4 x_5 \cov 
x_1 x_2 x_3 x_4 x_5  = {\hat 1}$.
\qed

The behavior of $\Lambda$ on intervals of rank two can be easily described.

\noindent{\bf Lemma 3.2.} 
{\it 
For every rank-two interval in a poset of shuffles $W_{M N}$, 
the labeling $\Lambda$ conforms to one of the following cases:
\hfil\break\indent
3.2.1)\ \  If a rank-two interval is isomorphic to $C_2 \times C_2$, 
then its two chains $c_1$ and $c_2$ have label sequences of the form
$\Lambda(c_1) = (l_1, l_2)$ and $\Lambda(c_2) = (l_2, l_1)$,
where $l_1$ and $l_2$ are distinct letters from ${\cal A} \cup {\cal X}$.
\hfil\break\indent
3.2.2)\ \  If a rank-two interval is isomorphic to $\Pi_3$, then its three
chains $\gamma_1$, $\gamma_2$ and $\gamma_3$ have label sequences of the 
form $\Lambda(\gamma_1) = (x_j, l)$, $\Lambda(\gamma_2) = (l, x_j)$, 
and $\Lambda(\gamma_3) = (x_j, x_j)$, for suitable letters 
$x_j \in {\cal X}$ and 
$l \in {\cal A} \cup ( {\cal X} - \{ x_j \} )$. 
}

{\bf Proof.} 
Each rank-two interval of a poset of shuffles has either 4 or 5 elements.
That is, each rank-two interval is isomorphic either to $C_2 \times C_2$
or to the lattice $\Pi_3$ of partitions of a 3-element set. 
Specifically,  an interval of rank 2 is of 
one of the following forms:
\hfil\break\indent 
i)\ \  $[ u a_m v a_n w, u v w]$ or $[u v w, u x_p v x_q w]$, 
for some shuffle words $u, v, w$;
\hfil\break
or $[ u a_m v w, u v x_p w ]$ or $[u v a_m w, u x_p v w ]$,
for some words $u, v, w$ with $v \neq \emptyset$.
By Lemma 2.2,  these intervals 
are isomorphic to $W_{1 0} \times W_{1 0}$, $W_{0 1} \times W_{0 1}$,
or $W_{1 0} \times W_{0 1}$, all of which are isomorphic to $C_2 \times C_2$.
\hfil\break\indent
ii)\ \ $[ u a_p v, u x_m v ]$ for some words $u, v$.
Such an interval is isomorphic to $W_{1 1} \simeq \Pi_3$. 

The definition of the chain labeling $\Lambda$
and the two possible structures of the intervals of rank 2
yield the two cases in the desired conclusion.
\qed

\noindent{\bf Proposition 3.3.}
{\it 
The labeling $\Lambda$ is an $S$-labeling on the poset of shuffles. 
}

{\bf Proof.} 
This follows immediately from Lemma 3.1 and Lemma 3.2.
\qed

\noindent{\bf Proposition 3.4.} 
{\it 
Consider the ordering  
$a_1 < a_2 < \cdots < a_M < x_1 < x_2 < \cdots < x_N$ 
on the union of the two alphabets.
Then the labeling $\Lambda$ is an $R^*$-labeling of the poset of shuffles. 
}

{\bf Proof.} 
Let 
${\hat 0} = w^0 \cov w^1 \cov \cdots \cov w^r < v$
be a chain in $W_{M N}$.
Let $B$ be the set of pairs of letters 
$x_{j_s}, a_{i_s}$ which occur as consecutive letters in $w^r$ 
and such that $a_{i_s} \not\in \{ v \}$. 
By the definition of $\Lambda$,  every covering along any $w^r$-$v$-chain 
where such an $a_{i_s}$ is removed
will be a covering of type ($xa$).  Thus, 
consider the $w^r$-$v$-chain obtained by first deleting, 
in increasing order of their indices, the letters in 
$({\cal A} - B) \cap ( \{ w^r \}  - \{ v \})$;  then 
inserting, also in increasing order of the indices, 
every $x_t \in \{ v \} - \{ w^r \}$, and deleting each letter $a_{i_s} \in B$ 
after the insertion of 
$x_t \in \{ v \} - \{ w^r \}$ if an only if $t < j_s$. 
The label sequence of this chain is clearly strictly increasing. 
Since any other label sequence with distinct entries 
is a permutation of the same set of labels $({\cal A}- B) \cup {\cal X}$,
this is the only strictly increasingly labeled $w^r$-$v$-chain. 
\qed

\noindent{\bf Remark 3.5.}
The reader familiar with the theory of shellable posets may be interested 
in the observation that 
{\it 
the labeling $\Lambda$ readily gives rise to a CL-labeling of $W_{M N}$
}
(in the sense of [BjWa], [BjGaSt]).
Indeed, the unique strictly increasing ${\hat 0}$-$u$-chain guaranteed 
by the preceding result can be taken as the ``root'' of each interval 
$[u,v]$, and the labeling $\Lambda^*$
defined by $\Lambda^*(u \cov v) = \left( \Lambda(u \cov v), M+N-\rho(u) \right)$
is a CL-labeling valued in 
$\left( ({\cal A} \cup {\cal X}) \times [M+N] \right)^{M+N}$, 
under lexicographic order on  
$({\cal A} \cup {\cal X}) \times [M+N]$.

\noindent{\bf Remark 3.6.}
{\it
The proof of Proposition 3.4 shows that the
$R^*$-labeling $\Lambda$ has a stronger property:
the unique increasingly labeled extension of a 
chain ${\hat 0} = w^0 \cov w^1 \cov \cdots \cov w^r < v$
depends only on $w^r$ and $v$.    
We will write $\gamma(w^r, v)$ to denote this chain.
}

The remainder of this section is devoted to enumerative consequences 
of the labeling $\Lambda$, yielding combinatorial proofs 
of results from [Gre1].
We begin with a bijective proof of the local rank-symmetry of the posets 
of shuffles (an inductive proof was given in Proposition 2.4). 
In particular, this is a bijective proof of the rank-symmetry of 
a poset of shuffles.   An alternative bijective proof of the rank-symmetry 
of $W_{M N}$  is implicit in the symmetric chain decomposition  
which appears in [Gre1].

\noindent{\bf Corollary 3.7.} 
{\it
For every two elements $u < w$ in a  poset of shuffles $W_{M N}$, 
there is a bijection between the elements of rank $\rho(u) + i$ and 
the elements of rank $\rho(w)-i$ in the interval $[u, w]$. 
}

{\bf Proof.}
Let $v \in [u, w]$ be an element of rank $\rho(u) + i$. 
Consider the maximal chain $c(u, v, w)$ formed by concatenating 
$\gamma({\hat 0}, u)$, $\gamma(u, v)$,  $\gamma(v, w)$, 
and $\gamma(w, {\hat 1})$.
Let $c'(u, v, w)$ be the unique maximal chain whose label sequence 
is the concatenation of 
$\Lambda(\gamma({\hat 0}, u))$, $\Lambda(\gamma(v, w))$, 
$\Lambda(\gamma(u, v))$,  
and $\Lambda(\gamma(w, {\hat 1}))$.
Define $\varphi(v)$ to be the element of rank $\rho(w) - i$ 
on the chain $c'(u,v,w)$.      
It is easy to see (from the definition and injectivity of $\Lambda$)
that $c'(u,v,w)$ contains the elements $u$ and $w$ 
and that $\varphi$ establishes a bijection between the 
rank-($\rho(u) + i$) and the rank-($\rho(w)-i$) elements in the 
interval $[u,w]$.
\qed

%
%

\noindent{\bf Corollary 3.8.}
{\it 
The number of elements of the poset $W_{M N}$ is equal to 
$\sum_{k \ge 0}{ {M \choose k} {N \choose k} 2^{M+N-2k}}$.
}

{\bf Proof.} 
We may count the increasingly labeled chains $\gamma({\hat 0}, w)$
since they are in bijection with the elements $w \in W_{M N}$.
For a prescribed number $k \ge 0$ of coverings of type $(xa)$, 
the set of labels along such a chain is determined by 
the choice of $k$ pairs from ${\cal A} \times {\cal X}$ 
for the coverings of type $(xa)$, and an arbitrary subset of the 
complement in ${\cal A} \cup {\cal X}$ of the letters chosen 
for the $k$ pairs.  
\qed

Lemma  3.1 yields readily the number of 
maximal chains in a poset of shuffles, giving a more direct derivation
of the formula  (6) due to C. Greene.

\noindent {\bf Corollary 3.9.} 
{\it
The number of maximal chains in the poset of shuffles $W_{M N}$ 
is 
$$
C_{M N} = \sum_{k \ge 0} { {M \choose k} {N \choose k} 
                           { {(M+N)!} \over {2^k} } }.
$$
}

{\bf Proof.} 
By Lemma 3.1, we can count the maximal chains in $W_{M N}$ by counting 
the possible label sequences  $\Lambda(c)$.  For each value  $k \ge 0$, 
the $k$th term in the sum gives the number of multipermutations $\Lambda(c)$
in which $k$ letters of $\cal X$ appear with multiplicity 2,
while $k$ of the letters of $\cal A$ do not occur in $\sigma$. 
\qed

{From} the $R^*$-labeling $\Lambda$ we can recover the formula (7) 
for the M\"obius function of a poset of shuffles,
which was obtained in [Gre1] using an EL-labeling as well as through 
an alternative computation. 

\noindent{\bf Corollary 3.10.} 
{\it 
The M\"obius function of the poset of shuffles $W_{M N}$ 
is given by}
 $$ \mu_{W_{M N}} ({\hat 0}, {\hat 1}) = (-1)^{M+N}{ {M+N} \choose M}.
 $$
{\bf Proof.} 
By the general theory of [BjW], 
the M\"obius function is $(-1)^{M+N}$ times the number of 
maximal chains to which the $R^*$-labeling $\Lambda$ assigns 
weakly decreasing label sequences. 
{From}  Lemma 3.1 it follows that such chains have 
label sequences of the form
$$
\Lambda(c) = ( x_{j_{N+k}}, x_{j_{N+k-1}}, \dots, x_{j_1}, 
   a_{i_{M-k}}, a_{i_{M-k-1}}, \dots, a_{i_1})
$$
 for some $0 \leq k \leq {\rm min} \{ M, N \}$, where 
$i_{M-k} > i_{M-k-1} > \cdots > {i_1}$
and 
$j_{N+k} \ge j_{N+k-1} \ge  \cdots  \ge {j_1}$ 
with $k$ non-consecutive equalities. 
Therefore, the decreasingly labeled maximal chains 
correspond bijectively to the selections of 
$M-k$ letters from $\cal A$ and $k$ letters from $\cal X$ for some $k$. 
It is an easy exercise to show that the number of 
such selections is ${ {M+N} \choose M}$ yielding the desired formula
for the M\"obius function.
\qed


\bigskip
\bigskip

\noindent{\bf 4.  A local action of the symmetric group}

We begin with two general results which imply that
the posets of shuffles have a local action of the symmetric group 
and establish the relation between 
the Frobenius characteristic of the corresponding character 
and the flag symmetric function of the poset.

\noindent{\bf Theorem 4.1.} 
{\it 
Suppose $P$ is a finite ranked poset of rank $n$, with $\hat 0$ and $\hat 1$. 
If $P$ has an $S$-labeling $\Lambda$, then the action of $S_n$ on
labels by permuting coordinates
induces a local (permutation) action on the maximal chains of $P$. 
}

{\bf Proof.} Let $\Lambda(c)=(\Lambda_1,\dots,\Lambda_n)$ be the label
of a maximal chain $c$, and let $1\leq i\leq n-1$. The adjacent
transposition $\sigma_i=(i,i+1)$ acts on $\Lambda(c)$ by interchanging
$\Lambda_i$ and $\Lambda_{i+1}$. By definition of $S$-labeling, there
is a unique maximal chain $c'$ such that $\Lambda(c')=\sigma_i\cdot
\Lambda(c)$. Since the $\sigma_i$'s generate $S_n$, we get an action
of $S_n$ on the set of labels of maximal chains, and hence on ${\cal
M}(P)$. Moreover, this action is local by the definition of an
$S$-labeling. \qed

\noindent{\bf Observation 4.2.1.} 
Suppose $S_n$ acts on the maximal chains of a labeled poset $P$ of rank $n$
by permuting the coordinates of the labels.
Then each orbit of maximal chains consists of the chains 
labeled by the permutations of a multiset, 
and the Frobenius characteristic of the $S_n$-action 
is $\sum_{\nu \vdash n}{N(\nu) h_{\nu}}$, where $N(\nu)$
denotes the number of orbits of maximal chains which are labeled 
by the permutations of a multiset of type $\nu$ 
(i.e., $\nu_1, \nu_2, \dots$ are the multiplicities of the 
distinct elements of the multiset). 

\noindent{\bf Observation 4.2.2.}
A special case is when the maximal chains in each orbit 
form a subposet isomorphic to a product of chains, 
$C_{\nu_1 + 1} \times C_{\nu_2 + 1} \times \cdots$.
It is not hard to show that this
is the case for the posets of shuffles and the action discussed here,
as well as for the lattice of noncrossing partitions discussed in [St3].
Thus, in addition to admitting a partition of 
the elements into boolean lattices (as shown in [Gre1] for poset 
of shuffles and in [SiU] for the noncrossing partition lattice), 
these posets also admit a partition of their maximal chains into 
products of chains.  
Figure 4 shows the  orbits of maximal chains in $W_{2 1}$.
In general, in $W_{M N}$, each orbit of maximal chains is isomorphic to 
a product of chains of the form $C_3^k \times C_2^{M+N-2k}$.

%
%
%
%

        \midinsert
        \epsfysize=3in
       \centerline{ \epsfbox{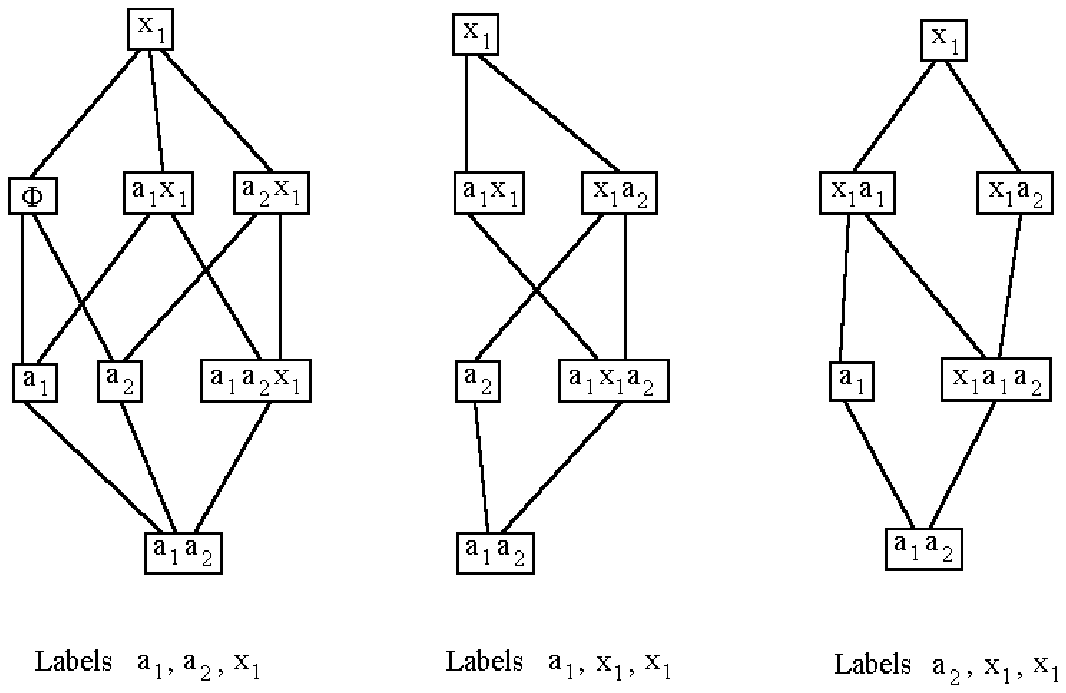} }
        \vskip -3ex
     \centerline{{\bf Figure 4.} The orbits of maximal chains in $W_{2
        1}$.} 
        \endinsert

Observations 4.2.1 and 4.2.2 are generalized by the following result.

\noindent{\bf Theorem 4.3.}  
{\it Suppose that $P$ is a ranked poset (with $\hat{0}$ and $\hat{1}$)
of rank $n$ with a local $S_n$-action.}
\hfil\break\indent
(a) {\it Let $c\in{\cal M}(P)$. Then the stabilizer stab$(c)$ of $c$
is a Young subgroup $S_{B_1}\times S_{B_2}\times \cdots\times S_{B_k}$
of $S_n$, where 
$\pi = \{B_1,B_2,\dots,B_k\}$ is a partition of $n$.}
\hfil\break\indent
(b) {\it If $\psi$ is the character of the $S_n$-action, then}
ch$(\psi)$ {\it is} $h$-positive, {\it i.e.,} ch$(\psi) =
\sum_{\nu\vdash n} a_\nu h_\nu$, {\it where} $a_\nu\geq 0$.
\hfil\break\indent
(c) {\it If} ch$(\psi)=h_\nu$ {\it for some $\nu\vdash n$, then
$P\simeq C_{\nu_1+1}\times C_{\nu_2+1}\times\cdots$.}

{\bf Proof.} (a) Let $\theta\in{\rm stab}(c)$, and let $i$ be the least
element of $[n]$ for which $\theta^{-1}(i)=j>i$. We claim that $(i,j) \in
{\rm stab}(c)$. Let $(i_1,i_2,\dots,i_r)$ be the cycle of $\theta$
containing $i$, where $i_1=i$ and $i_r=j$. Let $\tau_1,\tau_2,\dots,\tau_s$ be
the remaining cycles of $\theta$. Then (multiplying right-to-left)

\halign{\hfil#&#&#\hfil\cr
\ \hskip .1 true in   $\theta$& 
               $ = $ &$\tau_s\cdots \tau_2\tau_1(i_3,i_2)(i_4,i_3)
                        \cdots (i_r,i_{r-1})(i_1,i_r) $\cr
& $ = $ &$\tau_s\cdots \tau_2\tau_1(i_3,i_2)
     (i_4,i_3)\cdots (i_r,i_{r-1})\sigma_{i_r-1}\sigma_{i_r-2}
     \cdots \sigma_{i_1+1}\sigma_{i_1}\sigma_{i_1+1}\cdots
      \sigma_{i_r-2}\sigma_{i_r-1}.$\ \hskip .25 true in (11)\cr
}

\noindent Note that only one factor of the last product above moves
$i_1$, namely, $\sigma_{i_1}$. Let $t$ be the element of $c$ of rank
$i$. It follows from the definition of local action that $t$ is also
an element of the chain $c'=\sigma_{i_2}\cdots \sigma_{i_r-2}
\sigma_{i_r-1}\cdot c$. Let $s$ be the element of the chain $c''=
\sigma_{i_1}\sigma_{i_2}\cdots \sigma_{i_r-1}\cdot c$ of rank
$i$. Then again by definition of local action, $s$ is an element of
the chain $\theta \cdot c$ (since the factors to the left of $\sigma_{i_1}$
in (11) can be written as products of $\sigma_p$'s with $p>i$). Since
$\theta \cdot c=c$, we have $s=t$. Thus $\sigma_{i_1}\cdot c' =c'$, so we
can remove the factor $\sigma_{i_1}$ from the product (11) and still
get a permutation $\theta'\in {\rm stab}(c)$. But

\halign{\hfil #&#&#\hfil\cr
\ \hskip 1.75 true in $\theta'$ & $ = $ &$\tau_s\cdots \tau_2\tau_1 (i_3,i_2)
     \cdots (i_{r-1},i_{r-2})(i_r,i_{r-1})$\cr
& $ = $ &$\theta (i_1,i_r). $\cr
}

\noindent
Hence $(i_1,i_r)\in {\rm stab}(c)$, as claimed.

It follows by induction on $i$ (as defined above) that if for any
$a,b\in [n]$ we have $\theta(a)=b$, then $(a,b)\in {\rm stab}(c)$. From
this it is clear that stab$(c)$ is a Young  subgroup of $S_n$.

(b) By (a), the $S_n$-action on ${\cal M}(P)$, when restricted to an
orbit ${\cal O}\in {\cal M}(P)/S_n$, is equivalent to the action of
$S_n$ on the set $S_n/S_\nu$ of cosets of some Young subgroup $S_\nu =
S_{\nu_1}\times S_{\nu_2}\times\cdots$, where $\nu=\nu_{{\cal O}}
\vdash n$. If $\psi^\nu$ is the character of this action of $S_n$ on
$S_n/S_\nu$, then ch$(\psi^\nu)=h_\nu$. Hence
  $$ {\rm ch}(\psi) = \sum_{{\cal O}\in {\cal M}(P)/S_n} 
     h_{\nu_{{\cal O}}}. $$

(c) Let $M$ be the multiset $\{ 1^{\nu_1},2^{\nu_2},\dots\}$. The
action of $S_n$ on the set $S_M$ of permutations of $M$ obtained by
permuting coordinates is equivalent to the natural action of $S_n$ on
$S_n/S_\nu$. Hence by (a) there is an $S_n$-equivariant bijection
$\rho: {\cal M}(P)\rightarrow S_M$. Let $t\in P$, say rank$(t)=k$, and
let $c$ be a maximal chain of $P$ containing $t$. Let $\rho(c)=a=
a_1a_2\cdots a_n\in S_M$. Let $b=b_1b_2\cdots b_n\in S_M$ have the
property that $\{a_1,a_2,\dots,a_k\}=\{b_1,b_2,\dots,b_k\}$ (as
multisets), so also $\{
a_{k+1},\dots,a_n\}=\{b_{k+1},\dots,b_n\}$. Since $a$ can be
transformed to $b$ by adjacent transpositions all different from
$\sigma_k$, it follows from the definition of local action that the
chain $c'\in {\cal M}(P)$ satisfying $\rho(c')=b$ contains $t$. Hence
for any submultiset $N\subseteq M$, we can define $t_N$ to be the {\it
unique} element of $P$ for which there exists $c\in {\cal M}(P)$
containing $t_N$ and such that $N$ is equal to the first $k = {\rm
rank}(t_N)$ elements of $\rho(c)$. We thus have a well-defined
surjection $\tau: B_M\rightarrow P$, $\tau(N) = t_N$, where $B_M$ is
the lattice of submultisets of $M$ ordered by inclusion. Since $B_M
\simeq C_{\nu_1+1}\times C_{\nu_2+1}\times\cdots$, it suffices to show
that $\tau$ is a poset isomorphism.

By construction, $\tau$ is order-preserving (i.e., $N\subseteq N'
\Rightarrow \tau(N)\leq \tau(N')$), and the induced map $\tau: {\cal
M}(B_M) \rightarrow {\cal M}(P)$ is injective. Since $\#{\cal M}(B_M)
=\#{\cal M}(P) = n!/\nu_1!\,\nu_2!\cdots$, it follows that $\tau:{\cal
M}(B_M) \rightarrow {\cal M}(P)$ is a bijection. Suppose that $N,N'\in
B_M$ with $N\neq N'$ and $\tau(N)=\tau(N')$.  Let $c_1$ be a maximal
chain of the interval $[\emptyset,N]$ of $B_M$, and $c_2$ a maximal
chain of $[N',M]$. Then $\tau(c_1\cup c_2)$ is a maximal chain of $P$
not belonging to $\tau({\cal M}(B_M))$, contradicting the surjectivity
of $\tau: {\cal M}(B_M) \rightarrow {\cal M}(P)$. Thus $\tau$ is
injective on $B_M$.  Since $B_M$ and $P$ have the same number of
maximal chains and $\tau$ is injective, it is easy to see that $\tau$
must be an isomorphism. \qed


\noindent{\bf Theorem 4.4.} 
{\it 
Suppose $P$ is a flag-symmetric poset of rank $n$
with flag symmetric function $F_P$ 
and having an $S$-labeling $\Lambda$.  Let $\psi$ be the character 
of the action of $S_n$ on ${\bf C}{\cal M}(P)$ 
induced from the $S_n$-action on labels.  
\hfil\break\indent
(a)  If $\Lambda$ is an $RS$-labeling then $F_P = {\rm ch}(\psi) = h_{\nu}$
for some partition $\nu$ of $n$, 
and  $P \simeq C_{\nu_1 + 1} \times C_{\nu_2 + 1} \times \cdots$. 
\hfil\break\indent
(b)  If $\Lambda$ is an $R^*S$-labeling, then ${\rm ch}(\psi) = \omega
F_P$, where $\omega$ is the standard involution on symmetric functions
[M, p.21].
Hence $F_P$ is an $e$-positive symmetric function.
}

{\bf Proof.} (a)
Let $\gamma=(\gamma_1,\dots,\gamma_\ell)\in{\bf P}^\ell$, with
$\gamma_{1} +\cdots +\gamma_\ell=n$. For a finite multiset $M$, let
${\cal U}_\gamma(M)$ denote the collection of all sequences $\pi=(M_1,
\dots,M_\ell)$ of (nonempty) multisets $M_i$ such that
$\#M_i=\gamma_i$ and $\bigcup M_i=M$. Let 
  $$ {\cal U}_\gamma(\Lambda) = \bigcup_M {\cal U}_\gamma(M), $$
where $M$ ranges over all distinct multisets $M=\{\Lambda_1,\dots,
\Lambda_n\}$ of entries of labels $\Lambda(c)=(\Lambda_1,\dots,
\Lambda_n)$ of maximal chains of $P$. For each $\pi=(M_1,\dots,
M_\ell)\in {\cal U}_\gamma(\Lambda)$, it follows from the definition
of $RS$-labeling that there is a unique chain $t(\pi)=(\hat{0}=t_0 <
t_1<\dots <t_\ell=\hat{1})$ of $P$ with the following properties: 
\hfil\break\indent
(i)
$\rho(t_i) -\rho(t_{i-1})=\gamma_i$ for $1\leq i\leq \ell$, and 
\hfil\break\indent
(ii)
if $c$ is the unique completion of $t$ to a maximal chain of $P$ whose
label $\Lambda(c)=(\Lambda_1,\dots,\Lambda_\ell)$ satisfies
  $$ \hskip 1.3pc 
     \Lambda_1\leq\cdots\leq\Lambda_{\gamma_1},\ 
     \Lambda_{\gamma_1+1}\leq\cdots\leq\Lambda_{\gamma_1+\gamma_2},
     \dots,\ \Lambda_{\gamma_1+\cdots+\gamma_{\ell-1}+1}\leq\cdots
     \leq \Lambda_n, \hskip 2.5 pc \quad (12) $$
then $M_i=\{\Lambda_{\gamma_1+\cdots+\gamma_{i-1}+1},\dots,
\Lambda_{\gamma_1+\cdots \gamma_i}\}$. 

The map $\pi\mapsto t(\pi)$ is a bijection from ${\cal
U}_\gamma(\Lambda)$ to the set of all chains of $P$ whose elements
have ranks $0,\ \gamma_1,\ \gamma_1+\gamma_2,\dots,\ \gamma_1+\cdots +
\gamma_{\ell-1},\ n$. Hence
  $$ \hskip 5.5 pc
     F_P = \sum_{\lambda\vdash n} \alpha_P(\lambda)m_\lambda
     = \sum_{\lambda\vdash n} \#{\cal U}_\lambda(\Lambda)\cdot
     m_\lambda. \hskip 9 pc \qquad (13) $$
Let $\nu={\rm type}(M)$, i.e., $\nu\vdash n$ and the part
multiplicities of $M$ (in weakly decreasing order) are $\nu_1,\ \nu_2,
\dots$. It is well-known (equivalent to [M,(6.7)(ii)]) that
  $$ \hskip 9.2pc
     \sum_{\lambda\vdash n} \#{\cal U}_\lambda(M) \cdot m_\lambda
     = h_\nu. \hskip 10.5 pc \qquad (14) $$
Let ${\cal L}(\Lambda)$ be the collection of all multisets $M=\{
\Lambda_1,\dots,\Lambda_n\}$ of entries of a maximal chain label of
$P$. Summing (14) over all $M\in {\cal L}_\lambda$ and comparing with
(13) gives
  $$ F_P = \sum_{M\in{\cal L}(\Lambda)} h_{{\rm type}(M)}. $$
Since the action of $S_n$ by permuting coordinates of permutations of
multiset of type $\nu$ has Frobenius characteristic $h_\nu$, we get
$F_P = {\rm ch}(\psi)$. 

Since there is a {\it unique} weakly increasing maximal chain of $P$
from $\hat{0}$ to $\hat{1}$ (equivalently, since
$\alpha_P(\emptyset)=1$), we get $F_P=h_\nu$ for some $\nu\vdash n$.

It now follows from Theorem~4.3(c) that $P \simeq C_{\nu_1 + 1}
\times C_{\nu_2 + 1} \times \cdots$.  

(b) The argument is parallel to (a), except that the inequalities
$\leq$ of (12) become strict inequalities $<$. Hence ${\cal
U}_\gamma(M)$ is replaced by the collection ${\cal V}_\gamma(M)$ of
sequences $\pi = (M_1,\dots,M_\ell)$ of {\it sets}, rather than
multisets, so (14) becomes
  $$ \sum_{\lambda\vdash n} \#{\cal V}_\lambda(M) \cdot m_\lambda
     = e_\nu. $$
Since $\omega e_\nu=h_\nu$, we get $F_P=\omega({\rm ch}(\psi))$. \qed

The expression (10) for the flag symmetric function 
$F_{W_{M N}}$, follows now immediately from the general Theorem 4.4(b).

\noindent{\bf Corollary 4.5.} 
{\it 
The action of $S_{M+N}$ on the label sequences of the maximal chains 
in the poset of shuffles $W_{M N}$ induces a local action on the poset. 
The Frobenius characteristic for the 
character $\psi$ of the corresponding representation of $S_{M+N}$ 
 and the flag symmetric function of $W_{M N}$ are related by 
$$
{\rm ch}(\psi) = 
\sum_{k\ge 0}{ {M \choose k} {N \choose k} h_2^k(x) h_1^{M+N-2k}(x)} 
= \omega F_{W_{M N}}(x).
$$
}

In the remainder of this section we make comments regarding the 
preceding results and discuss other possible directions for generalizations.

\noindent{\bf Remark 4.6.} 
Theorems 4.1 and 4.4 apply to posets which are products of chains 
and also to the lattice of noncrossing partitions.
The corresponding conclusions are established directly in [St2] and [St3].

\noindent{\bf Remark 4.7.} 
The power series $F_P(x)$ may be viewed in a broader context.
For a function $\varphi$ in the incidence algebra 
(see, e.g., [St1, section 3.6])  
of a ranked poset $P$,  define 
$$\alpha_P(\varphi, S) = 
     \sum_{{\hat 0} = t_0 < t_1 < \cdots < t_k < {\hat 1}}
          { \varphi({\hat 0}, t_1) \varphi(t_1, t_2) 
                                                \cdots \varphi(t_k, {\hat 1})
          },
$$
where the sum ranges over the chains in $P$ whose rank-support is 
the set $S \subseteq [n-1]$, and $n$ is the rank of $P$.
Now define
$$
F_P(\varphi, x) \ \colon = 
 \sum_{{S \subseteq [n-1]} \atop{ S = \{ s_1 < s_2 < \cdots < s_k \} } } 
      { \alpha_P(\varphi, S)   \cdot 
       {\sum_{1 \leq i_1 < i_2 < \cdots < i_{k+1} } 
            {x_{i_1}^{s_1} x_{i_2}^{s_2 - s_1} \cdots x_{i_{k+1}}^{n - s_k} }
        }
      }.
$$
Note that $F_P(\varphi, x)$ is a quasisymmetric function, 
homogeneous of degree $n$.
In particular,  if $\varphi = \zeta$, the zeta function of $P$
(i.e.,  $\zeta(u, v) = 1$ if $u \leq v$, and $\zeta(u,v) = 0$ otherwise),
we recover $\alpha_P(\zeta, S ) = \alpha_P(S)$ and 
 $F_P(\zeta, x)$ is the function $F_P(x)$ of (2).
We intend to pursue this generalization elsewhere, mentioning here 
only one result -- the next proposition.   
We note that the same result holds for an arbitrary invertible 
element $\varphi$ from the incidence algebra of $P$ and its 
inverse.  Here we present a direct proof for the 
special case $\varphi = \zeta$ and $\varphi^{-1} = \mu$ 
which is the instance occurring in the context of this paper.
 
\noindent{\bf Proposition 4.7.1.}
{\it 
Let $P$ be a ranked poset of rank $n$, having elements $\hat 0$ and $\hat 1$.
If $\zeta$ and $\mu$ are, as usual, the zeta function 
and the M\"obius function of $P$, then 
$$
F_P(\mu, x) = (-1)^n  \omega F_P(\zeta, x),
$$
where $\omega$ is the involution on quasisymmetric functions 
defined by $\omega L_{S,n}(x) = L_{ {\overline{S}}, n}(x)$, 
with ${\overline{S}}$ denoting the complement of $S$ in $[n-1]$.
}
 
{\bf Proof.} 
Using (3) and then (1), we have 
$$
\omega F_P(\zeta, x) = \omega \sum_{S \subseteq [n-1]}{\beta_P(S) L_{S,n}(x)} 
                     = \sum_{S \subseteq [n-1]}
                           { \sum_{T \subseteq S}{(-1)^{|S - T|} \alpha_P(T)}
                             L_{{\overline{S}}, n}(x) }.
$$ 
Using Hall's theorem (e.g., [St1, Prop. 3.8.5]), 
the sum over $T$ evaluates to $(-1)^{|S|-1} \mu_{P_S}({\hat 0}, {\hat 1})$.
Next, using Baclawski's theorem for the M\"obius function of 
a subposet (see [Gre2, formula (7.2)]), the sum over $T$ can be 
expressed as 
$$
(-1)^{|S|-1} \sum_{k \ge 0}
            {\sum_{
                   {{\hat 0} < t_1 <  t_2 < \cdots < t_k < {\hat 1}}
                    \atop {{\rm in \ } P_{\overline{S}} }
                   }
             {(-1)^k \mu({\hat 0}, t_1) \mu(t_1, t_2) \cdots \mu(t_k,{\hat1})}}.
$$
(When $k=0$, the chain ${\hat 0} < {\hat 1}$ gives the term 
$\mu({\hat 0},{\hat 1})$.)
By grouping the chains in $P_{\overline{S}}$ according to their 
rank-support, $U \subseteq {\overline{S}}$, we obtain 
$$
\omega F_P(\zeta, x) = (-1)^n 
                      \sum_{U \subseteq [n-1]}
                       { \alpha_P(\mu, U) \cdot 
                         \sum_{ {\overline{S}} \supseteq U} 
                          {(-1)^{|{\overline{S}} - U|} L_{{\overline{S}}, n}(x)}
                       }.
$$
Finally, by the definition of $L_{V,n}(x)$ and an inclusion-exclusion argument,
the inner sum over ${\overline{S}}$ is equal to 
$\sum_{1 \leq i_1 < i_2 < \cdots < i_j}{x_{i_1}^{u_1} x_{i_2}^{u_2 - u_1} 
                \cdots x_{i_{j+1}}^{n - u_j}}$,
where $U = \{ u_1 < u_2 < \cdots < u_j \}$. 
Thus, $\omega F_P(\zeta, x) = (-1)^n F_P(\mu, x)$ as claimed. 
\qed

Since $\omega$ restricted to symmetric functions agrees with the standard 
involution $\omega$, the expressions for the characteristic ch$(\psi)$ 
from Theorem 4.4 can be restated as 
(a) \ ${\rm ch}(\psi) = F_P(\zeta, x) = h_{\nu}$ when $P$ is $RS$-labeled,
and (b) \ ${\rm ch}(\psi) = (-1)^n F_P(\mu, x)$ when $P$ is $R^*S$-labeled.

\noindent{\bf Remark 4.8.} 
{\it
The $S_{M+N}$-action of Corollary 4.5  is a permutation action on
the maximal chains of $W_{M N}$, for which each orbit consists of those 
maximal chains whose label sequence $\Lambda(c)$ 
is a permutation of the same multiset of letters.
Thus, the explanation of formula (6) for the maximal chains of $W_{M N}$ 
provided in the proof of Corollary 3.9 amounts to 
counting the maximal chains according 
to their orbit, and grouping  the orbits according to the 
type $2^k 1^{M+N-2k}$ of the multiset of the chain labels. 
}

%
%

      \midinsert
        \epsfysize=3in
       \centerline{ \epsfbox{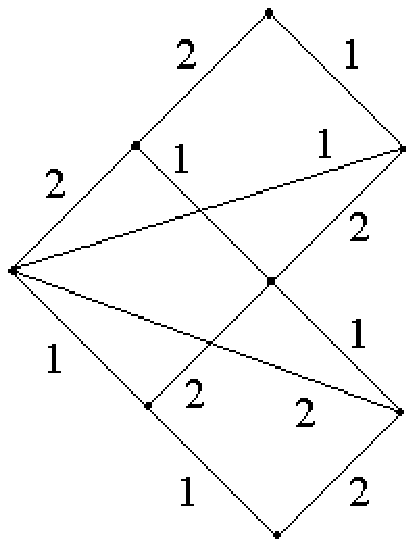} }
        \vskip -3ex
\centerline{{\bf Figure 5.} A poset with an $S$-labeling}
\centerline{and an orbit of maximal chains which is not a product of chains.}
        \endinsert

\noindent{\bf Remark 4.9.} 
{\it 
An $S$-labeling does not ensure that each orbit of maximal
chains is isomorphic to a product of chains. 
}
The poset shown in Figure 5, suggested to us 
by H. Barcelo, has a local action induced from 
the action of $S_4$ on the labels, but the 
orbit of chains labeled by the multiset $1 1 2 2$ is not a product of 
chains.  In fact, no product of chains other than 
the trivial one, $C_5$, occurs as a subposet of rank four in this poset.

\noindent{\bf Remark 4.10.} 
{\it 
The converse of Theorem 4.1 
does not hold. That is, a chain labeling such that the action of $S_n$
on labels induces a local action is not necessarily an $S$-labeling. 
}
The poset $P$ shown in Figure~6(a) 
has a labeling of its maximal chains which, although not injective, 
gives rise to an $S_3$  local action on the maximal chains of $P$. 
The orbits are four copies of $C_2 \times C_3$, 
each labeled in the standard way with the multiset $ a, b, b $.
On the other hand, the noninjective labeling of Figure~6(b) 
does not produce an $S_3$ action on the maximal chains 
(e.g., $\sigma_1 \sigma_2 \sigma_1({\hat 0} \cov A \cov B \cov {\hat 1})
\neq \sigma_2 \sigma_1 \sigma_2 ({\hat 0} \cov A \cov B \cov {\hat 1})$).

%
%
%
%

        \midinsert
        \epsfysize=3in
       \centerline{ \epsfbox{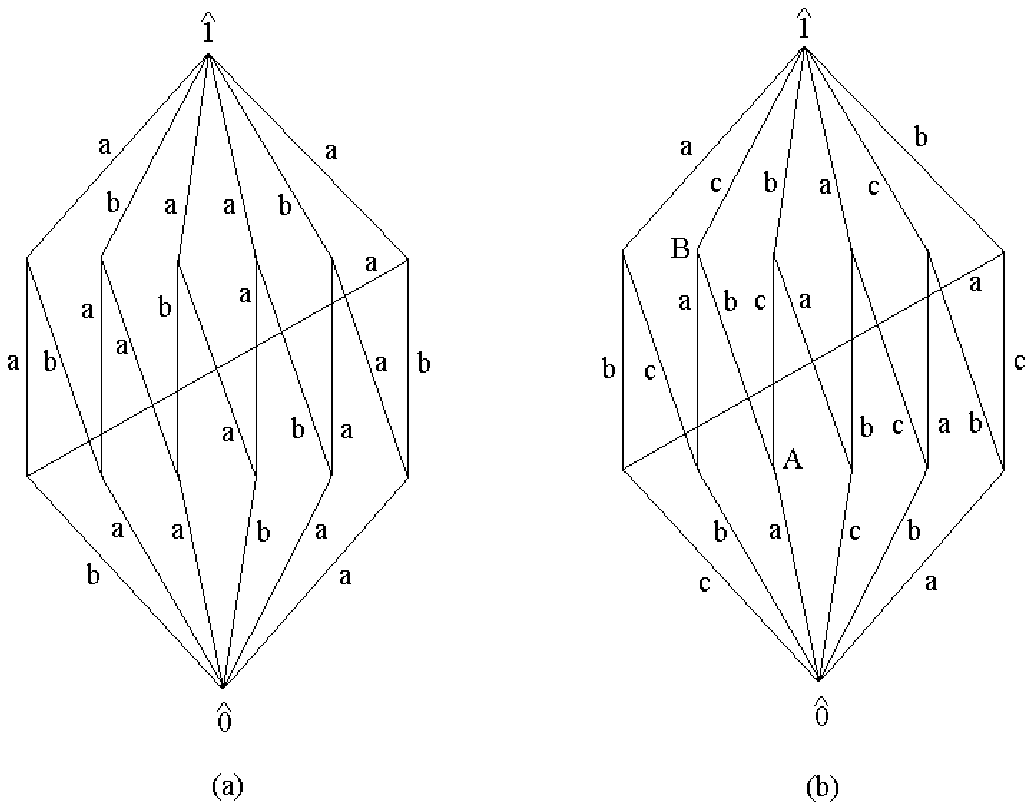} }
\centerline{{\bf Figure 6.}  Non-injective labelings.}
\centerline{a) The action on labels induces a local action.}
\centerline{b) The action on labels does not induce a local action.}
        \endinsert

Finally, we give a local condition which characterizes labeled posets 
with a local action induced from the action on labels.
This, of course, can be seen to apply to the earlier examples. 

\noindent{\bf Theorem 4.11.} 
{\it 
Let $P$ be a finite ranked poset of rank $n$, having a $\hat 0$ and $\hat 1$,
and with a labeling $\Lambda$ of its  maximal chains.
For a maximal chain 
$c = ({\hat 0} = t_0 \cov t_1 \cov \cdots \cov t_n = {\hat 1})$ and 
a rank $i \in \{ 1, 2, \dots, n-2 \}$, 
let $\tau \colon =  ({\hat 0} = t_0 \cov t_1 \cov \cdots \cov t_{i-1})$ 
and 
let $\theta \colon =  (t_{i+2} \cov t_{i+3} \cov \cdots \cov t_n = {\hat 1})$.
A local action is induced from the $S_n$-action on labels if and only if 
$\Lambda$ satisfies the following condition for 
every maximal chain $c$ and every value of $i \in \{ 1, 2, \dots, n-2 \}$: 
\hfil\break\indent
The length-three chains $\delta$ from $t_{i-1}$ to $t_{i+2}$
with labels induced by restricting $\Lambda(\tau \delta \theta)$ 
can be partitioned so that  
\hfil\break\indent
(a)  each class is isomorphic to 
$C_2 \times C_2 \times C_2$ or $C_3 \times C_2$ or $C_4$, and 
\hfil\break\indent
(b) the labeling in each class coincides with the standard labeling 
of a product of chains by join-irreducibles.
}

{\bf Proof.} 
The conditions (a) and (b) on $\Lambda$ imply readily the Coxeter
relations for $S_n$, showing that the local action is well-defined.
Conversely, within each orbit of chains $\delta$, the local action
fixes the chains labeled $a a a$, so these form classes isomorphic to
$C_4$; a chain $\delta$ labeled with $a b a$ is mapped under the local
action to chains with the same $\tau$ and $\theta$ and label sequences
$a a b$ and $b a a$, structured as a copy of $C_2 \times C_3$ and
forming a class as claimed; similarly, a chain $\delta$ labeled as $a
b c$ is mapped by the subgroup generated by $\sigma_i$ and $\sigma_{i+1}$ to six
chains forming a copy of $C_2 \times C_2 \times C_2$, with labels as
claimed.
\qed


\bigskip
\bigskip


\noindent{\bf 5. Multiplicative functions on the poset of shuffles}

Consider now the poset $W_{\infty\infty}$ whose elements are the
shuffles of finite words using the lower alphabet ${\cal A}_{\infty} =
\{ a_1, a_2, \dots \}$ and the upper alphabet ${\cal X}_{\infty} = \{
x_1, x_2, \dots \}$.  The comparability relation is as in the case of
finite alphabets.  A {\it multiplicative function} on
$W_{\infty\infty}$ is a function $f$ defined on the intervals in
$W_{\infty \infty}$ for which $f_{00}=1$ and which has the property
that if $[ u, v ] \caniso \prod_{i,j}{ W^{c_{ij}}_{ij} }$ then $f(u,v)
= \prod_{ij}{f^{c_{ij}}_{ij}}$, where we write $f_{ij}$ for the value
of $f$ on an interval canonically isomorphic to $W_{ij}$ (see Remark
2.3).

Let $f$ and $g$ be two multiplicative functions on $W_{\infty \infty}$,
and let 
\baselineskip = 2 pc
\halign{\hfil#&#&\hfil#&#&#\hfil\cr
\hskip 1.4 true in
$F$ & \ $=$ \  &$F(x,y)$ & \ $=$ \ &$\sum_{i,j\geq 0} {f_{ij}\,x^iy^j}$\cr
$G$ & \ $=$ \  &$G(x,y)$ & \ $=$ \ &$\sum_{i,j\geq 0} {g_{ij}\,x^iy^j}$\cr
$F * G$ & \ $=$ \ &$(F * G)(x,y)$ & \ $=$ \ &$\sum_{i,j\geq 0} 
                                              {(f * g)_{ij}\,x^iy^j},$ \cr
}
\baselineskip = 1.3 pc
\noindent
where $f * g$ denotes convolution in the incidence algebra
$I(W_{\infty  \infty})$. The main result of this section, Theorem 5.2, 
shows how to express $F * G$
in terms of $F$ and $G$, and hence ``determines'' the monoid of
multiplicative functions on $W_{\infty \infty}$.

We begin by establishing an expression for the number of elements 
$w \in W_{M N}$ of a given {\it type} $( (a_{ij}), (b_{ij}))$,
that is, such that 
   $[ {\hat{0}}, w] \caniso \prod_{i,j}{W_{ij}^{a_{ij}}}$
   and 
   $[ w, {\hat{1}}] \caniso \prod_{i,j}{W_{ij}^{b_{ij}}}$.
Note that the canonical isomorphism (Remark 2.3) 
implies that 
$1+\sum ib_{ij}=\sum a_{ij}$ and $1+\sum ja_{ij}=\sum b_{ij}$,
and we can recover from the type of a word $w$ the values  
$m \colon = \# (\{w\} \cap {\cal A})  
          = \sum_{i,j}{i \left(a_{ij} + b_{ij}\right)}$ 
and 
$n \colon = \# (\{w\} \cap {\cal X})  
          = \sum_{i,j}{j \left(a_{ij}+b_{ij}\right)}. $ 
Also, the type of $w \in W_{M N}$ determines $M$ and $N$,
so the enumeration of shuffle words by type can be done in $W_{\infty \infty}$.

\noindent{\bf Proposition 5.1.} 
{\it
   Let $(a_{ij})$ and $(b_{ij})$ be nonnegative integers
   such that 
   $$ \epsilon \ \colon = \sum_{{i,j} \atop {j \neq 0}}{a_{ij}} - 
              \sum_{{i,j} \atop {i \neq 0}}{b_{ij}} \  \in \{ 1, 0, -1
              \}. $$
   Set 
   \hfil\break\indent
   $$ m = \sum_{i,j}{i \left(a_{ij} + b_{ij}\right)}, \quad 
      n = \sum_{i,j}{j \left(a_{ij}+b_{ij}\right)}, $$ 
  and 
$r = \sum_{{i,j}\atop{j\ne 0}}{a_{ij}}$, and \ \  
$s = \sum_{{i,j}\atop{i\ne 0}}{b_{ij}}$. 
 
   Then the number  of elements $w \in W_{\infty \infty}$  whose type is 
 $( (a_{ij}), (b_{ij}))$ is given by 
$$
   \cases{
   (2 - \epsilon^2) \displaystyle
    { { { {m+1} \choose {\dots, a_{ij}, \dots } } 
        { {n+1} \choose {\dots, b_{ij}, \dots } }
      }
    \over
      { { {m+1} \choose r }
        { {n+1} \choose s } 
      } 
    } & if $w\neq\emptyset$ (i.e., $r+s > 0$)\cr
    1& if $w=\emptyset$ (i.e., $r=s=0$).\cr
    }$$
}

\noindent{\bf Proof.} 
Each $w \in W_{\infty \infty} - \{ \emptyset \}$ 
is of the form $U L U L \cdots U L$, 
or $L U L U \cdots L U$, or $L U \cdots L U L$, 
or $U L \cdots U L U$, 
where each $U$ is a nonempty factor whose letters are from the
upper alphabet, and each $L$ is a nonempty factor whose letters are
from the lower alphabet.  If the type of $w$ is $( (a_{ij}),
(b_{ij}))$ then, for each $j \ge 1$, the number of $U$-factors
of length $j$ is $\sum_{i}{a_{ij}}$ and, for each $i \ge 1$, the
number of $L$-factors of length $i$ in $w$ is equal to
$\sum_{j}{b_{ij}}$.  The alternation of nonempty $L$- and $U$-factors imposes
the condition $\epsilon \in \{ 1, 0, -1 \}$ appearing in the
hypothesis.  To construct a word $w$ of the prescribed type, we begin
by deciding the length of each $U$- and $L$-factor.  The number of
possibilities is the number of (multi)permutations of the nonzero
lengths, so that $U$- and $L$-factors alternate:
$$
\hskip 0.2 true in
{ {r \choose 
  {a_{01}, a_{02}, \dots , a_{11}, a_{12}, \dots }}} 
{ {s \choose 
  {b_{10}, b_{20}, \dots , b_{11}, b_{21}, \dots }}} 
  (1 + \chi(\epsilon  = 0)),
\hskip .7 true in (15)
$$
where, following A. Garsia [Ga] (see also [Kn]), 
if $p$ is a proposition then we write $\chi(p)=1$ if $p$
is true, and $\chi(p)=0$ if $p$ is false. 
To complete the construction of $w$, we need to 
choose the location of the $W_{i0}$'s and $W_{0j}$'s required by 
entries $a_{i0}$ and $b_{0j}$ in the type of $w$. 
A factor $W_{i0}$ in the canonical product for $[{\hat 0}, w]$
must arise between two successive lower alphabet  
letters of $w$, or in front of the first $L$-factor if $w$ begins 
with an $L$-factor, 
or after the last $L$-factor if $w$ ends with an $L$-factor. 
Therefore such a factor $W_{i0}$ occurs in one of  
$m - s + 1 - \epsilon$  positions.
Similarly, a factor $W_{0j}$ in the canonical product for 
$[ w, {\hat 1}]$ can arise from any of $n - r + 1 + \epsilon$  positions 
(between two successive letters from the upper alphabet, 
 in front of the first $U$-factor if $w$ begins 
with a $U$-factor, or after the last $U$-factor if $w$ ends with a
$U$-factor).  
%
%
In conclusion,   the word $w$ can be completed in 
$$
{ {m - s + 1 - \epsilon} 
  \choose {a_{10}, a_{20}, \dots,  m - s + 1 - \epsilon - \sum_{i}{a_{i0}}}}
{ {n -r + 1 + \epsilon} 
  \choose {b_{01}, b_{02}, \dots,  n - r + 1 + \epsilon - 
    \sum_{j}{b_{0j}}}}
$$
ways.

The remainder of the proof is a calculation. 
After multiplying (15) by the preceding expression, 
the relations 
$r-s = \epsilon$,  
$1+\sum ib_{ij}=\sum a_{ij}$, and $1+\sum ja_{ij}=\sum b_{ij}$
allow some simplifications.
For example, $m -s + 1 - \epsilon - \sum_i{a_{i0}} = 
m + 1 - r - \sum_i{a_{i0}} = 
\sum_i{i b_{ij}} + 1 - (r + \sum_i{a_{i0}}) = 0$.
Similarly, $n - r + 1 + \epsilon - \sum_j{b_{0j}} = 0$. 

The first case in the 
conclusion of the Proposition now follows from a simple manipulation
with binomial coefficients. The case $w=\emptyset$ is trivial, so the
proof is complete.  
\qed


\noindent{\bf Theorem 5.2.} 
{\it
\baselineskip = 2 pc
Let $F_0=F(x,0)$, $G_0=G(0,y)$, and
\halign{\hfil#&#&#\hfil\cr
\hskip 1.9 true in$\tilde{F}(x,y)$ &\ $ =$ \ &$F(x,G_0y)$\cr
                  $\tilde{G}(x,y)$ &\ $ =$ \ &$G(F_0x,y)$.\cr
}
\baselineskip = 1.3 pc
Then
$$
\hskip 1.85 true in
{1 \over {F * G}} = {1 \over {\tilde{F}G_0}}+{1 \over {F_0 \tilde{G}}}
                   -{1 \over {F_0G_0}}. 
\hskip 1.95 true in  (16)
$$
}

{\bf Proof.} For fixed $r,s,m,n \geq 0$ write
  $$ Q_{r,s,m,n} = \sum {
                          {\prod{f_{ij}^{a_{ij}}g_{ij}^{b_{ij}} 
                                 x^{\sum{i(a_{ij} + b_{ij})}}
                                 y^{\sum{j(a_{ij} + b_{ij})}}}
                          }
                         \over
                          {( \prod{a_{ij}!} )( \prod{b_{ij}!} )}
                          }, 
  $$
where the sum ranges over all $a_{ij}$ and $b_{ij}$ satisfying
 $$ \sum_{j\neq 0}a_{ij}=r,\quad
                          \sum_{i\neq 0}b_{ij}=s $$
  $$ \sum ib_{ij}=m,\quad \sum ja_{ij} =n $$
  $$ m+1=\sum a_{ij}, \quad n+1=\sum b_{ij}. $$
By Proposition 5.1, the convolution $F * G$ is given by 
$$ 
F * G=(F * G)_{-1} + (2(F * G)_0-F_0G_0)+(F * G)_1, 
$$ 
where
\baselineskip = 2 pc
\halign{\hfil#&#&#\hfil\cr
\hskip .6 true in 
$(F * G)_{-1}$ &\ $ = $\
&$\displaystyle\sum_{k,m,n}{k!\, (k+1)!\, (m-k)!\, (n+1-k)!\,
Q_{k+1,k,m,n}} $\cr 
$(F * G)_0$ &\ $ = $\ &$\displaystyle\sum_{k,m,n}{ k!^2\, (m+1-k)!\,
(n+1-k)!\, Q_{k,k,m,n}}$\cr 
$(F * G)_1$ &\ $ = $\ &$\displaystyle\sum_{k,m,n}{(k+1)!\, k!\,
(m+1-k)!\, (n-k)!\, Q_{k,k+1,m,n}}$.\cr
}
\baselineskip = 1.3 pc
If $M$ is a monomial then write $[M]Q$ for the coefficient of $M$ in
the power series $Q$. 
We first consider $(F * G)_0$. 
We have

\baselineskip = 2.3 pc
\halign{#&#&#\hfil\cr
$(F * G)_0$ &\ $ = $\ &$\displaystyle\sum_{k,m,n}
                   {k!^2\, (m +1-k)!\, (n+1-k)!\,
                   [q^{m+1}r^{n+1}s^mt^nu^kv^k]}$\cr 
            & & \ $\cdot \displaystyle\sum_{a_{ij},b_{ij}}
   {{{\prod_{j\neq 0} {\left( q\,t^j uf_{ij}x^iy^j\right)^{a_{ij}}}}
     {\prod_i {\left( q\, f_{i0}x^i\right)^{a_{i0}}}}
     {\prod_{i\neq 0} {\left(rs^ivg_{ij}x^iy^j\right)^{b_{ij}}}}
     {\prod_j {\left(rg_{0j}y^j\right)^{b_{0j}}}}
    } 
    \over
   {\left( \prod{a_{ij}!} \right) \left( \prod{b_{ij}!} \right)
    } 
   }$\cr\cr
             &\ $ = $\ &$\displaystyle \sum_{k,m,n}
                         {k!^2\, (m+1-k)!\, (n+1-k)!
                   [q^{m+1}r^{n+1}s^mt^nu^kv^k]}$\cr 
             &         &$ \ \cdot \displaystyle \prod_{i,j}
                  { \left(
  \displaystyle \sum_{a_{ij}\geq 0}{ {\left( q\,t^j\,u^{\chi(j\neq
                   0)}f_{ij}\,x^iy^j\right)^{a_{ij}}} 
                       \over
                       {a_{ij}!}
                      }\right)
                  \cdot 
                  {\left(
  \displaystyle\sum_{a_{ij}\geq 0}{ {\left( rs^iv^{\chi(i\neq
                   0)}g_{ij}x^iy^j\right)^{b_{ij}}} 
                       \over
                       {b_{ij}!}
                      }\right) }  }$\cr\cr
             &\ $ = $\ &$\displaystyle\sum_{k,m,n} {k!^2\,
                   (m+1-k)!\,(n+1-k)!}$\cr 
             &         &$\ \cdot
  [q^{m+1}r^{n+1}s^mt^nu^kv^k] 
   \displaystyle \prod_{i,j}{\exp{\left( q\, t^j u^{\chi(j\neq
                   0)}f_{ij}\, x^iy^j 
                   + rs^iv^{\chi(i\neq 0)}g_{ij}x^iy^j\right)}}$\cr\cr
            &\ $ = $\ &$\displaystyle \sum_{k,m,n}
                   {k!^2\,(m+1-k)!\,(n+1-k)!}$\cr 
            &         &$ \ \cdot
  [q^{m+1}r^{n+1}s^mt^nu^kv^k] 
  \exp \displaystyle\sum_{i,j}{\left( q\,t^j u^{\chi(j\neq 0)}f_{ij}\,
                   x^iy^j 
                  + r\,s^iv^{\chi(i\neq 0)}g_{ij}\,x^iy^j\right)}$\cr\cr
            &\ $ = $\ &$\displaystyle \sum_{k,m,n}{
                   k!^2\,(m+1-k)!\,(n+1-k)!}$\cr 
            &         &$ \ \cdot
   [s^mt^nu^kv^k] { \displaystyle
                   {\left( \sum_{i,j}
                      {t^ju^{\chi(j\neq 0)}f_{ij}\,x^iy^j}\right)^{m+1}
                    }
                    \over 
                    {\displaystyle (m+1)!}
                 } 
                  \cdot
                  {
                   {\displaystyle\left( \sum_{i,j} 
                      {s^iv^{\chi(i\neq 0)}g_{ij}\,x^iy^j}\right)^{n+1}
                    }
                    \over
                    {\displaystyle (n+1)!}
                   }$\cr\cr
           &\ $ = $\ &$\displaystyle \sum_{k,m,n} 
                  {
                   {k!^2\,(m+1-k)!\,(n+1-k)!}
                  \over
                   {(m+1)!\,(n+1)!}
                 } 
                 [s^mt^n] {{m+1}\choose k} 
                 \left( \displaystyle\sum_i{f_{i0}\,x^i} \right)^{m+1-k}$\cr
           &          &$ \cdot 
                 \left( \displaystyle\sum_{{i,j}\atop {j\neq 0}}
                                {f_{ij}\,t^jx^iy^j}  \right)^k 
                {\displaystyle{{n+1}\choose k}} 
                 \left( \displaystyle\sum_i{g_{0j}\,y^j} \right)^{n+1-k}
                 \left( \displaystyle\sum_{{i,j}\atop{i\neq
                   0}}{g_{ij}\,s^ix^iy^j} 
                 \right)^k$\cr\cr
           &\ $ = $\ &$\displaystyle\sum_{k,m,n} 
                       {[s^mt^n] 
                        \left( \displaystyle\sum_i{f_{i0}\,x^i}
                   \right)^{m+1-k} 
                        \left( \displaystyle\sum_{{i,j}\atop{j\neq 0}}
                                   {f_{ij}\,t^jx^iy^j} \right)^k 
                        \left( \displaystyle\sum_i {g_{0j}\,y^j}
                   \right)^{n+1-k} 
                        \left( \displaystyle\sum_{{i,j}\atop{i\neq 0}} 
                                   {g_{ij}\,s^ix^iy^j} \right)^k} $\cr\cr
           &\ $ = $\ &$ \displaystyle \sum_k 
                       { \left(
                   \displaystyle\sum_i{f_{i0}\,x^i}\right)^{1-k}  
                         \left(
                   \displaystyle\sum_j{g_{0j}\,y^j}\right)^{1-k} 
                       }$\cr
           &         &$\ \cdot
                         \displaystyle\sum_{m,n}
                         {[s^mt^n]  
                          \left( \displaystyle\sum_{{i,j}\atop{j\neq 0}}
                                     {f_{ij}\, t^jx^iy^j} \right)^k  
                          \left( \displaystyle\sum_{{i,j}\atop{i\neq 0}}
                                     {g_{ij}\, s^ix^iy^j} \right)^k
                          \cdot 
                          \left( \displaystyle\sum_i{f_{i0}\, x^i} \right)^m
                          \left( \displaystyle\sum_j{g_{0j}\, y^j} \right)^n
                          } $\cr\cr
           &\ $ = $\ &$ \displaystyle\sum_k
                       {\left( \displaystyle\sum_i{f_{i0}\,x^i} \right)^{1-k}
                        \left( \displaystyle\sum_j{g_{0j}\,y^j} \right)^{1-k}
                        }$\cr 
           &         &$\ \cdot 
                       \left( \displaystyle\sum_{{i,j}\atop{j\neq 0}}
                                 {\left(
                   \displaystyle\sum_r{g_{0r}\,y^r}\right)^j  
                                   f_{ij},x^iy^j} 
                       \right)^k
                       \left( \displaystyle\sum_{{i,j}\atop{ i\neq 0}}
                                 {\left(
                   \displaystyle\sum_s{f_{s0}\,x^s}\right)^i  
                                   g_{ij}\,x^iy^j} 
                       \right)^k.$\cr
}
\baselineskip = 1.3 pc
Note that
  $$ \sum_r{g_{0r}\,y^r} =G_0 $$
and
  $$ \sum_{{i,j}\atop{j\neq 0}}{z^jf_{ij}\,x^iy^j} = F(x,yz)-F_0, $$
and similarly for $F_0$ and $G(xz,y)-G_0$. Hence
\baselineskip = 2 pc
\halign{\hfil#&#&#\hfil\cr
\hskip .25 true in $(F * G)_0$
    &\ $ = $\ &$F_0 G_0\cdot 
                \displaystyle\sum_k{F_0^{-k}G_0^{-k}
                       \left( \displaystyle\sum_{{i,j}\atop{j\neq 0}}
                                      {G_0^j\,f_{ij}\,x^iy^k} 
                       \right)^k
                       \left( \displaystyle\sum_{{i,j}\atop{i\neq 0}}
                                      {F_0^i\,g_{ij}\,x^iy^k}
                       \right)^k} $\cr
    &\ $ = $\ &$\displaystyle{ {F_0G_0}
                \over
                  {1- { {(\tilde{F}-F_0)(\tilde{G}-G_0)} \over {F_0G_0} }}
                }.$\cr
}
\baselineskip = 1.3 pc
Exactly analogous reasoning applies to $(F * G)_{-1}$ and $(F * G)_1$. 
For instance, $(F * G)_{-1}$ can be written
\baselineskip = 2.3 pc 
\halign{#&#&#\hfil\cr
$(F * G)_{-1}$ &\ $ = $\ &$\displaystyle\sum_{k,m,n}{k!\,(k+1)!\,
                              (m-k),(n+1-k)!\,     
                                        [q^{m+1}r^{n+1}s^mt^nu^{k+1}v^k]}$\cr
               &         &$ \ \cdot \displaystyle\sum_{a_{ij},b_{ij}}
                                     { 
                     {\prod_{j\neq 0}{\left( q\,t^j
                              u\,f_{ij}\,x^iy^j\right)^{a_{ij}}} 
                      \prod_i{\left( q\,f_{i0}\,x^i\right)^{a_{i0}}}
                      \prod_{i\neq
                              0}{\left(r\,s^iv\,g_{ij}\,x^iy^j\right)^{b_{ij}}} 
                      \prod_j{\left(r\,g_{0j}\,y^j\right)^{b_{0j}}} }
                     \over
                     {\left( \prod{a_{ij}!} \right)
                      \left( \prod{b_{ij}!} \right)}
                                    } .$\cr
}
\baselineskip = 1.3 pc
Reasoning as above yields
  $$ (F * G)_{-1} = { {(\tilde{F}-F_0)G_0}
                      \over
                      {1- { {(\tilde{F}-F_0)(\tilde{G}-G_0)}\over{F_0G_0}}}
                     }. $$
Similarly (or by symmetry), 
  $$ (F * G)_1 = {  {F_0(\tilde{G}-G_0)}
                 \over
                    {1- { {(\tilde{F}-F_0)(\tilde{G}-G_0)}\over{F_0G_0}}}
                 }. $$
{From} this we easily obtain
$$
{1 \over {F * G}}\  = \        {1 
                                \over
                                 {(F * G)_{-1} +(2(F * G)_0-F_0G_0)+(F * G)_1}
                                }
                    \  = \    {1 \over {\tilde{F}G_0}}
                               +{1 \over {F_0 \tilde{G}}}
                               -{1 \over {F_0G_0}}, 
$$
as desired.
\qed

\noindent{\bf Example 5.3.} In general it seems difficult to understand the
operation $F * G$. It is not even obvious from (16) that
$*$ is associative! A special case for which $F * G$ can be
explicitly evaluated is the following. Let $a_1, b_1, \dots, a_k, b_k$
be real numbers (or indeterminates). Then it is straightforward to
prove from Theorem~5.2 by induction on $k$ that
$$
{1 \over{(1-a_1x)(1-b_1y)}} * \cdots * 
{1 \over{(1-a_kx)(1-b_ky)}} =
\hskip 1 true in 
$$
$$
\hskip 1 true in
{1 \over {1-(\sum{a_i})x - (\sum{b_i})y
          + (\sum{(a_1+a_2+\cdots+a_i)b_i}) xy}
 }. \hskip 1.2 true in (17)
$$
For instance, if $\zeta$ denotes the zeta function of
$W_{\infty \infty}$ (whose value is 1 on every interval of
$W_{\infty \infty}$), then
$$
\sum_{i,j}{\zeta_{ij}x^iy^j} = 
  {1 \over {(1-x)(1-y)}}. 
$$
Hence the left-hand side of (17) becomes the generating
function for $\zeta^k$, whose value at $W_{ij}$ is the number
$Z_{ij}(k)$ of $k$-element multichains ${\hat 0}=t_0\leq t_1\leq \cdots
\leq t_k={\hat 1}$ in $W_{i j}$. Equivalently, $Z_{ij}(k)$ is the value
of the zeta polynomial of $W_{i j}$ at $k$ [St1].  We obtain
from (17) that
$$ 
\sum_{i,j}{Z_{ij}(k)x^iy^j} = 
{1 \over {1-kx-ky-{{k+1}\choose 2}xy}},
$$
a result of Greene [Gre1].

\noindent{\bf Remark 5.4.} 
{\it 
The identity (17) can be deduced purely combinatorially.
}
\hfil\break\indent
The left-hand side provides a refinement of the flag $f$-vector
in the sense that the coefficient of 
$a_1^{i_1} b_1^{j_1} a_2^{i_2} b_2^{j_2} \cdots x^m y^n$,
where $m = i_1 + i_ 2 + \cdots$ and $n = j_1 + j_2 + \cdots$, 
is the number of multichains ${\hat 0} = t_0 \leq t_1 \leq \cdots 
\leq t_{m+n-1} < t_{m+n} = {\hat 1}$ in $W_{m n}$
such that $t_{r-1}$ and $t_r$ differ by $i_r$ letters from 
the alphabet ${\cal A}$ and by $j_r$ letters from the alphabet ${\cal X}$. 
The right-hand side of (17) can be interpreted as the 
generating function for the language ${\cal L}$ over the alphabet 
$\{ a_1, a_2, \dots, b_1, b_2, \dots \}$ consisting of the words 
in which there is no occurrence of a letter $a_k$ immediately followed
by a letter $b_l$ with $k \leq l$. 
(This follows directly through a simple sign-reversing-involution 
argument or from the general theory of Cartier and Foata [CaFo].)
The number of such words 
formed with the multiset of letters 
$a_1^{i_1} b_1^{j_1} a_2^{i_2} b_2^{j_2}\cdots $,
where $m = i_1 + i_ 2 + \cdots$ and $n = j_1 + j_2 + \cdots$, 
is the coefficient of 
$a_1^{i_1} b_1^{j_1} a_2^{i_2} b_2^{j_2} \cdots x^m y^n$ 
on the right-hand side of (17).   
For example, the coefficient of $a_1^2 b_1 b_2 x^2 y^2$ is $2$, 
accounting for the words $b_1 b_2 a_1 a_1$ and $b_2 b_1 a_1 a_1$ in $\cal L$.

Now each word $w \in {\cal L}$ consisting of  $m$ $a_i$'s and $n$ $b_i$'s
determines a shuffle-word $s(w)$ in $W_{m n}$ 
by placing the alphabet ${\cal A}$ 
in the positions of the $a_i$'s and the alphabet ${\cal X}$ in the positions 
of the $b_i$'s.  For example, $w = a_2 b_3 b_3 a_1 a_3 b_5 b_1 b_2 a_3$ 
gives rise to the shuffle word 
$s(w) = a_1 x_1 x_2 a_2 a_3 x_3 x_4 x_5 a_4$ in $W_{4 5}$. 
Of course, $w \to s$ is not a one-to-one map. 

Given a word $w \in {\cal L}$, we construct a unique multichain 
$t(w)$ in $W_{m n}$ by starting with $t_0 = {\hat 0}$ and 
obtaining $t_{r+1}$ from $t_r$ by removing the letters appearing in $s(w)$ 
in those positions where $a_r$ occurs in $w$, and inserting the 
letters appearing in $s(w)$ in those positions where $b_r$ occurs 
in $s(w)$.   For instance, 
the word $w$ from the preceding example yields 
$t(w) = ({\hat 0}= t_0 < t_1 < t_2 < t_3 = t_4 < t_5 = {\hat 1})$,
where 
$t_0 = a_1 a_2 a_3 a_4$,
$t_1 = a_1 a_3 x_4 a_4$,
$t_2 = a_3 x_4 x_5 a_4$,
$t_3 = t_4 = x_1 x_2 x_4 x_5$,
$t_5 = x_1 x_2 x_3 x_4 x_5$. 

Recall from Section 3
that there is a unique order in which the $i_r$ deletions and $j_r$ 
insertions can be performed such that $t_{r+1}$ is reached from $t_r$ 
via the $\Lambda$-increasing chain $\gamma(t_r, t_{r+1})$.
The condition defining the 
words $w \in {\cal L}$ ensures that each insertion of a letter from ${\cal X}$ 
is made in the leftmost position permitted by the shuffle word $s(w)$,
thus not multiply-counting chains which involve covering relations of 
the type $(xa)$. 
{From} these observations it follows that 
$w \to t(w)$ is a bijection, completing a combinatorial proof of (17).

\bigskip
\bigskip

\noindent{\bf References}

\item{[BjGaSt]} A. Bj\"orner, A. Garsia and R. Stanley, 
{\it An introduction to Cohen-Macaulay partially ordered sets,} 
in ``Ordered Sets,'' NATO Adv. Study Inst. Ser. C: Math. Phys. Sci. 
(I. Rival, Ed.), Reidel, Dordrecht/Boston, 1982, pp. 583-615.
\item{[BjW]} A. Bj\"orner and M. Wachs, 
{\it On lexicographically shellable posets,} Trans. Amer. Math. Soc. 
{\bf 277} (1983), 323-341. 
\item{[CaFo]} P. Cartier and D. Foata, ``Probl\`emes Combinatoires
de Commutation et R\'earrangements,'' 
Lecture Notes in Math., no. 85, Springer-Verlag, Berlin/Heidelberg/New York, 
1969.
\item{[E]} R. Ehrenborg, {\it Posets and Hopf algebras}, Advances in
Math.\ {\bf 119} (1996), 1--25.
\item{[FaSc]} J. Farley and S. Schmidt, 
{\it Posets that locally resemble
distributive lattices: an extension of Stanley's theorem,} in preparation.
\item{[Ga]} A. Garsia, {\it On the `maj' and `inv' q-analogues
of Eulerian polynomials,} Lin.\ and Multilin. Alg. {\bf 8} (1979) 21-34. 
\item{[Gra]} D. Grabiner,  
{\it Posets in which every interval is a product of chains, and
natural local action of the symmetric group,} \hfil\break manuscript,
1997.  
\item{[Gre1]} C. Greene, {\it Posets of shuffles,} 
J. Combin. Theory Ser. A {\bf 47} (1988), 191-206.
\item{[Gre2]} C. Greene, {\it The M\"obius function of a partially 
ordered set,} in 
``Ordered Sets,'' NATO Adv. Study Inst. Ser. C: Math. Phys. Sci. 
(I. Rival, Ed.), Reidel, Dordrecht/Boston, 1982, pp. 555-581.
\item{[He]} P. Hersh, personal communication.
\item{[Kn]}  D. Knuth, {\it Two notes on notation,} 
Amer. Math. Monthly {\bf 99} (1992), 403-422.  
\item{[M]} I. G. Macdonald, `` Symmetric Functions and Hall Polynomials ,'' 
second edition, Oxford University Press, Oxford, 1995.   
\item{[Re]} C. Reutenauer, ``Free Lie Algebras,'' Oxford University
Press, New York, 1993.
\item{[SiU]} R. Simion and D. Ullman, {\it On the structure of 
the lattice of non-crossing partitions,} Discr. Math. {\bf 98} (1991), 193-206.
\item{[St1]} R. Stanley, ``Enumerative Combinatorics,'' volume 1, 
Wardsworth \& Brooks/Cole, Monterey, CA, 1986; second printing, Cambridge
University Press, Cambridge/New York, 1997.
\item{[St2]} R. Stanley, {\it Flag-symmetric and locally rank-symmetric 
partially ordered sets,} Electronic J. Combinatorics {\bf 3}, R6 (1996), 22 pp.
\item{[St3]} R. Stanley, {\it Parking functions and noncrossing partitions}, 
Electronic J. Combinatorics {\bf 4}, R20 (1997), 17 pp.
\item{[V]} A. M. Vershik, {\it Local stationary algebras,} 
Amer. Math. Soc. Transl. (2) {\bf 148} (1991), 1-13; translated from 
Proc. First Siberian Winter School ``Algebra and Analysis'' 
(Kemerovo, 1988). 

\vfil\eject
\bye